\documentclass[12pt]{amsart}
\usepackage{amsmath,amssymb,amsfonts,amsthm,amsopn}
\usepackage{graphics}

\newcommand{\C}{\mathbb{C}}

\newcommand{\modsp}{modulation space}

\newtheorem{tm}{Theorem}[section]
\newtheorem{lemma}[tm]{Lemma}
\newtheorem{prop}[tm]{Proposition}

\newtheorem{Remark}[tm]{Remark}

\newtheorem{theorem}{Theorem}[section]
\newtheorem{corollary}[theorem]{Corollary}

\newtheorem{proposition}[theorem]{Proposition}

\newcommand{\beqa}{\begin{eqnarray*}}
\newcommand{\eeqa}{\end{eqnarray*}}

\newcommand{\field}[1]{\mathbb{#1}}
\newcommand{\bR}{\field{R}}        
 \def\cF{\mathcal{F}}              
 
 \def\cD{\mathcal{D}}
 
 \def\cB{\mathcal{B}}

 \def\cC{\mathcal{C}}

\def\a{\aleph}

\def\rd{\bR^d}

\def\lrd{L^2(\rd)}

\def\intrd{\int_{\rd}}

\def\R{\right)}

\def\<{\left<}
\def\>{\right>}

\def\inv{^{-1}}

\def\mv1{M_v^1}

\def\o{\xi}
\def\a{\alpha}

\def\ZZ{\mathbb{Z}}

\def\R{\mathbb{R}}
\def\Ren{\mathbb{R}^d}
\def\Renn{\mathbb{R}^{2d}}

\def\sch{\mathcal{S}}

\def\Fur{\mathcal{F}}

\def\Sn2{S_{2}(L^{2}(\Ren))}
\def\S1{S_{1}(L^{2}(\Ren))}
\def\sig00{\sigma_{0,0}}

\newcommand{\fpq}{W(\mathcal{F}L^p,L^q)}

\begin{document}

\title[]{Some new Strichartz estimates
for the Schr\"odinger
equation}

\author{Elena Cordero and Fabio Nicola}
\address{Department of Mathematics,  University of Torino,
Via Carlo Alberto 10, 10123
Torino, Italy}
\address{Dipartimento di Matematica, Politecnico di
Torino, Corso Duca degli
Abruzzi 24, 10129 Torino,
Italy}
\email{elena.cordero@unito.it}
\email{fabio.nicola@polito.it}



\date{}

\begin{abstract} We deal with
fixed-time and Strichartz  estimates for the Schr\"o\-din\-ger
propagator as an operator on Wiener
 amalgam spaces. We discuss
 the sharpness of the known estimates
and we provide some new
estimates which generalize
the classical ones.
  As an application, we present
  a result on the
  wellposedness of
 the linear Schr\"odinger equation with
 a rough time dependent potential.
\end{abstract}

\maketitle

\section{Introduction}
Consider the Cauchy problem for the
Schr\"odinger equation
\begin{equation}\label{cp}
\begin{cases}
i\partial_t u+\Delta u=0\\
u(0,x)=u_0(x),
\end{cases}
\end{equation}
with $x\in\R^d$, $d\geq1$.
 Several estimates have
been obtained for the
solution
$u(t,x)=\left(e^{it\Delta}u_0\right)(x)$
of \eqref{cp} in terms of
Lebesgue spaces, with
important applications to
wellposedness and scattering
theory for
 nonlinear Schr\"odinger equations,
possibly with potentials
\cite{Dancona05,GinibreVelo92,Kato70,keel,strichartz,tao2,tao,Yajima87}.
Among them, we recall the important fixed-time estimates
\begin{equation}\label{est3}
\|e^{it\Delta}u_0\|_{L^{r}(\bR^d)}\lesssim
|t|^{-d\left(\frac{1}{2}-\frac{1}{r}\right)}\|u_0\|_{L^{r^\prime}(\bR^d)},\quad
2\leq r\leq\infty,
\end{equation}
as well as the homogeneous
Strichartz estimates: for
$q\geq2$, $r\geq 2$, with
$2/q+d/r=d/2$,
$(q,r,d)\not=(2,\infty,2)$,
i.e., for $(q,r)$ {\it
Schr\"odinger admissible},
\begin{equation}\label{S1}
\|e^{it\Delta} u_0\|_{L^q_t
L^r_x}\lesssim \|u_0\|_{L^2_x},
\end{equation}
where, as usual, $\displaystyle
\|F\|_{L^q_tL^r_x}=\left(\int\|F(t,\cdot)\|^q_{L^r_x}\,dt\right)^{1/q}.$
As a matter of fact, these estimates express  a gain of local
$x$-regularity of the solution $u(t,\cdot)$, and a decay of its
$L^r_x$ norm,
 both in some $L^q_t$-averaged
 sense.\par
Recently several authors
(\cite{benyi,
benyi2,cordero,cordero2,baoxiang3,baoxiang2,baoxiang})
have turned their attention
to fixed-time and space-time
estimates of the
Schr\"odinger propagator
between spaces widely  used
in Time-Frequency Analysis
\cite{grochenig}, known as
modulation spaces and Wiener
amalgam spaces. These spaces
were first introduced by H.
Feichtinger in
\cite{feichtinger80} (now
available in textbooks
\cite{grochenig}) and have
recently  become very popular
in the framework of signal
analysis. Loosely speaking,
given a Banach space $B$,
like  $B=L^p$ or $B=\Fur
L^p$, the Wiener amalgam
space $W(B,L^q)$ consists of
functions which are locally
in $B$ but display an
$L^q$-decay at infinity. In
particular, $W(L^p,L^p)=L^p$,
$W(L^{p_1},L^{q_1})\hookrightarrow
W(L^{p_2},L^{q_2})$ if
$p_1\geq p_2$, $q_1\leq q_2$,
and $W(\Fur
L^{p_1},L^{q_1})\hookrightarrow
W(\Fur L^{p_2},L^{q_2})$ if
$p_1\leq p_2$, $q_1\leq q_2$.
\par As explained in
\cite{cordero,cordero2}, one
of the main motivations for
considering estimates in
these spaces is that they
control the local regularity
of a function and its decay
at infinity separately.
Hence, they highlight and
distinguish between the local
properties and the global
behaviour of the solution
$u(t,x)=e^{it\Delta}u_0$ and
therefore, as far as
 Strichartz estimates concern, they are
natural candidates to perform
an analysis of the solution
which is finer than the
classical one.
 Actually,
Wiener amalgam spaces have already  appeared  as technical tool in
PDE's. In particular, the space $W(L^p,L^q)$ coincides with the
space $X^{q,p}_0$ introduced in \cite{tao3}.\par

Some Stricharz estimates in this environment are in
\cite{cordero,cordero2}, where   the following fixed-time
estimates are presented: for $2\leq r\leq \infty$,
\begin{equation}\label{est2}
\|e^{it\Delta}u_0\|_{W(\Fur
L^{r^\prime}, L^r)}\lesssim
\left(|t|^{-2}+|t|^{-1}\right)^{d\left(\frac{1}{2}-\frac{1}{r}\right)}\|u_0\|_{W(\Fur
L^r,L^{r^\prime})}.
\end{equation}
The related homogeneous Strichartz estimates, obtained by combining \eqref{est2} with orthogonality arguments, read
\begin{equation}\label{str}
\|e^{it\Delta} u_0\|_{W(L^{{q}/{2}},L^q)_t W(\Fur
L^{r^\prime},L^r)_x}\lesssim \|u_0\|_{L^2_x},
\end{equation}
for $4<q\leq\infty$, $2\leq
r\leq\infty$, with $2/q+d/r =d/2$. When
$q=4$ the same estimate holds with the
Lorentz space $L^{r',2}$ in place of
$L^{r'}$.  Dual homogeneous  and
retarded estimates hold as well.\par
For comparison, the classical estimates
\eqref{S1} can be rephrased in terms of
Wiener amalgam spaces as follows:
\begin{equation}\label{S1W}
\|e^{it\Delta}
u_0\|_{W(L^q,L^q)_t
W(L^r,L^r)_x}\lesssim
\|u_0\|_{L^2_x}.
\end{equation}
Thereby, the new estimates \eqref{str}
contain the following insight: the
classical estimates \eqref{S1W} can be
modified by (conveniently) moving local
regularity from the time variable to
the space variable. Indeed, $\Fur
L^{r'}\subset L^r$ if $r\geq2$, but the
bound in \eqref{est2} is worse  than
the one in  \eqref{est3}, as $t\to0$;
consequently, in the estimate
\eqref{str} we average locally in
 time by the  $L^{q/2}$ norm, which is rougher
 than the $L^q$ norm in \eqref{S1} or, equivalently, in \eqref{S1W}.\par
 In the
present paper we perform the converse
approach, by showing that it is
possible to move local regularity in
\eqref{S1} from the space variable to
the time variable. As a result, we
obtain new estimates involving the
Wiener amalgam spaces $W(L^p,L^q)$ that
generalize \eqref{S1}.  This requires
some preliminary steps.\par First, we
show that
 \eqref{est2} can be enlarged to   more general pairs $(s',r)$, rather
 than only    conjugate-exponent pairs $(r',r)$, (see Theorem \eqref{tfix}).
In particular, if we choose $s'=2$, the related  fixed-time
estimates read
\begin{equation*}
\|e^{it\Delta}u_0\|_{W(
L^{2}, L^r)}\lesssim (1+
t^2)^{-\frac{d}2\left(\frac{1}{2}
-\frac{1}{r}\right)}\|u_0\|_{W(L^{2},
L^{r^\prime})},\quad 2\leq
r\leq\infty.
\end{equation*}
Using  techniques similar to \cite{cordero,keel}, the related Strichartz estimates are achieved in Theorem \ref{prima2}. Finally, the complex
interpolation with the classical estimates \eqref{S1} yields our
main result:
\begin{theorem}\label{prima2i}
Let $1\leq q_1,r_1\leq\infty$, $2\leq
q_2,r_2\leq\infty$ such that $r_1\leq
r_2$,
\begin{equation}\label{pri1}\frac{2}{q_1}+
\frac{d}{r_1}\geq\frac{d}{2},
\end{equation}
\begin{equation}\label{pri2}\frac{2}{q_2}+
\frac{d}{r_2}\leq\frac{d}{2},
\end{equation}
$(r_1,d)\not=(\infty,2)$,
$(r_2,d)\not=(\infty,2)$ and,
if $d\geq 3$, $r_1\leq
2d/(d-2)$. The same for
$\tilde{q}_1,\tilde{q}_2,\tilde{r}_1,\tilde{r}_2$.
Then, we have the homogeneous
Strichartz estimates
\begin{equation}\label{hom2i}\|e^{it\Delta}
u_0\|_{W(L^{q_1},L^{q_2})_t
W(
L^{r_1},L^{r_2})_x}\lesssim
\|u_0\|_{L^2_x},
\end{equation}
the dual homogeneous
Strichartz estimates
\begin{equation}\label{dh2i}
\|\int e^{-is\Delta}
F(s)\,ds\|_{L^2}\lesssim
\|F\|_{W(L^{\tilde{q}_1^\prime},L^{\tilde{q}_2^\prime})_t
W(
L^{\tilde{r}_1^\prime},L^{\tilde{r}_2^\prime})_x},
\end{equation}
and the retarded Strichartz
estimates
\begin{equation}\label{ret2i}
\|\int_{s<t} e^{i(t-s)\Delta}
F(s)\,ds\|_{W(L^{q_1},L^{q_2})_t
W( L^{r_1},L^{r_2})_x}
\lesssim\|F\|_{W(L^{\tilde{q}_1^\prime},L^{\tilde{q}_2^\prime})_t
W(
L^{\tilde{r}_1^\prime},L^{\tilde{r}_2^\prime})_x}.
\end{equation}
\end{theorem}
Figure 1 illustrates the
range of exponents for the
homogeneous estimates when
$d\geq3$. Notice that, if
$q_1\leq q_2$, these estimates
follow immediately from \eqref{S1W} and the
inclusion relations of Wiener
amalgam spaces recalled
above. So, the issue consists
in the cases $q_1>q_2$.

  \vspace{1.2cm}

 \begin{center}
           \includegraphics{figSharp.1}
            \\
           $ $
\end{center}
           \begin{center}{Figure 1:
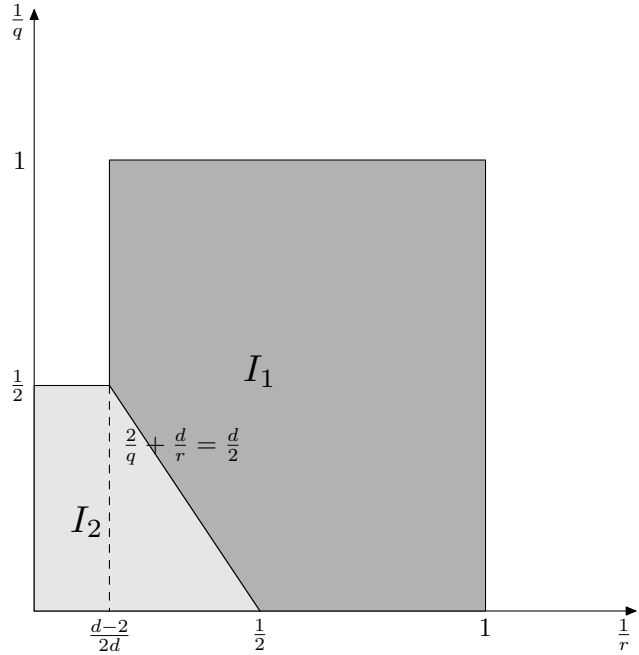
           When $d\geq3$, \eqref{hom2i} holds for all pairs $(1/q_1,1/r_1)\in
I_1$, $(1/q_2,1/r_2)\in I_2$,
with $1/r_2\leq 1/r_1$.}
           \end{center}
  \vspace{1.2cm}

Since there are no relations between
the pairs $(q_1,r_1)$ and $(q_2,r_2)$
other than $r_1\leq r_2$, these
estimates tell us, in a sense, that the
analysis of the local regularity of the
Schr\"odinger propagator is quite
independent of its decay at
infinity.\par

We then discuss the sharpness
of the results above, as well
as those in \cite{cordero}.
Indeed, in Section
\ref{necessarie} we first
focus on the fixed-time
estimates, proving that the
range $r\geq2$ in
\eqref{est2} is sharp, and
the same for the decay
$t^{-d\left(\frac{1}{2}-\frac{1}{r}\right)}$
at infinity, and the bound
$t^{-2d\left(\frac{1}{2}-\frac{1}{r}\right)}$,
when $t\to 0$. Then, we
present the sharpness of the
Strichartz estimates
\eqref{str}, except for the
threshold $q\geq4$, which
seems quite hard to obtain.
Next, we turn our attention
to the new estimates in
Theorem \ref{prima2i} and
show that, for $d\geq3$, all
the constraints on the range
of exponents in Theorem
\ref{prima2i} are necessary,
except for $r_1\leq r_2$,
$r_1\leq 2d/(d-2)$, which
remain an open problem.
However, we prove the weaker
result (Proposition
\ref{pd1}):\par
\medskip \emph{Assume
$r_1>r_2$ and $t\not=0$. Then
the propagator $e^{it\Delta}$
does not map
$W(L^{r'_1},L^{r'_2})$
continuously into
$W(L^{r_1},L^{r_2})$.}
\medskip

\noindent This shows that estimates
\eqref{hom2i} for exponents $r_1>r_2$,
if true, cannot be obtained from
fixed-time estimates and orthogonality
arguments.
\par The arguments here employed for the necessary conditions
differ from the classical setting of
Lebesgue spaces, where necessary
 conditions  are usually
obtained by general scaling considerations
(see, for example,
\cite[Exercises 2.35, 3.42]{tao},  and \cite{nicola} for the
interpretation in terms of Gaussian curvature of the
characteristic manifold). In our framework
 this method does not
work directly. Indeed, the known bounds
for the norm of the dilation operator
$f(x)\longmapsto f(\lambda x)$ between
 Wiener amalgam
spaces (\cite{sugimototomita,Toft04}), yield constraints which are
weaker than the desired ones.

Our idea is to  consider families of
rescaled Gaussians as initial data, for
which the action of the operator
$e^{it\Delta}$ and the involved norms
can be computed explicitly.\par In the
last section, we present an application
to the linear Schr\"odinger equation
with time-dependent potential (see
\cite{Dancona05}
 and the references therein for the existent
literature). Our result
 extends that in \cite{cordero} to  the dimension $d\geq 1$ (instead
of $d>1$) and to potentials $V(t,x)$ in Wiener amalgam spaces
rather than the classical $L^p$ spaces. Precisely, we prove the
wellposedness in $L^2$ of the Cauchy problem
\begin{equation}
\begin{cases}
i\partial_t u+\Delta u=V(t,x)u,\quad t\in [0,T]=I_T, \,\,x\in\rd,\\
u(0,x)=u_0(x),
\end{cases}
\end{equation}
for the class of potentials
\begin{equation}
V\in L^\alpha(I_T;W(\cF
L^{p^\prime},
L^p)_x),\quad\frac1{\a}+\frac{d}{p}\leq1,\
1\leq\alpha<\infty,\ {d}<
p\leq\infty.
\end{equation}
This result seems of interest
especially in dimension $d=1$, where,
if we choose $1< p <2$,   $V(t,x)$ is
allowed to be locally in rough spaces
of temperate distributions $\Fur
L^{p'}$, with respect to the space
variable $x$ (see Remark \ref{espot}).
\par
Estimates similar to those proved here should hold for other
dispersive equations, like the wave equation, too. Our plan is to
investigate these issues  in a subsequent
paper.\par\medskip\noindent \textbf{Notation.} We define
$|x|^2=x\cdot x$, for $x\in\Ren$, where $x\cdot y=xy$ is the
scalar product on $\Ren$. The space of smooth functions with
compact support is denoted by $\cC_0^\infty(\rd)$, the Schwartz
class is $\sch(\Ren)$, the space of tempered distributions
$\sch'(\Ren)$.    The Fourier transform is normalized to be ${\hat
  {f}}(\o)=\Fur f(\o)=\int
f(t)e^{-2\pi i t\o}dt$.
 Translation and modulation operators ({\it time and frequency shifts}) are defined, respectively, by
$$ T_xf(t)=f(t-x)\quad{\rm and}\quad M_{\o}f(t)= e^{2\pi i \o
 t}f(t).$$
We have the formulas
$(T_xf)\hat{} = M_{-x}{\hat
{f}}$, $(M_{\o}f)\hat{}
=T_{\o}{\hat {f}}$, and
$M_{\o}T_x=e^{2\pi i
x\o}T_xM_{\o}$. The notation
$A\lesssim B$ means $A\leq c
B$ for a suitable constant
$c>0$, whereas $A \asymp B$
means $c\inv A \leq B \leq c
A$, for some $c\geq 1$. The
symbol $B_1 \hookrightarrow
B_2$ denotes the continuous
embedding of the linear space
$B_1$ into $B_2$.
\section{Wiener amalgam spaces}
We briefly recall the
definition and the main
properties of Wiener amalgam
spaces. We refer to
\cite{feichtinger80,feichtinger83,
feichtinger90,fournier-stewart85,Fei98})
for details.\par
 Let $g \in
\cC_0^\infty$ be a test
function that satisfies
$\|g\|_{L^2}=1$. We will
refer to $g$ as a window
function. Let $B$ one of the
following Banach spaces:
$L^p, \cF L^p$, $1\leq p\leq
\infty$,
 $L^{p,q}$, $1<p<\infty$, $1\leq q\leq \infty$,
  possibly valued in a Banach space, or also
  spaces obtained from these by real or
  complex interpolation.
Let $C$ be one of the
following Banach spaces:
$L^p$, $1\leq p\leq\infty$,
or $L^{p,q}$,
 $1<p<\infty$, $1\leq q\leq \infty$, scalar-valued.
For any given function $f$
which is locally in $B$ (i.e.
$g f\in B$, $\forall
g\in\cC_0^\infty$), we set
$f_B(x)=\| fT_x g\|_B$.

The {\it Wiener amalgam
space} $W(B,C)$ with local
component $B$ and global
component  $C$ is defined as
the space of all functions
$f$ locally in $B$ such that
$f_B\in C$. Endowed with the
norm
$\|f\|_{W(B,C)}=\|f_B\|_C$,
$W(B,C)$ is a Banach space.
Moreover, different choices
of $g\in \cC_0^\infty$
generate the same space and
yield equivalent norms.

If  $B=\Fur L^1$ (the Fourier
algebra),  the space of
admissible windows for the
Wiener amalgam spaces $W(\Fur
L^1,C)$ can be enlarged to
the so-called Feichtinger
algebra $W(\Fur L^1,L^1)$.
Recall  that the Schwartz
class $\sch$
  is dense in $W(\Fur L^1,L^1)$.\par
We use the following
definition of mixed Wiener
amalgam norms. Given a
measurable function $F$ of
the two variables $(t,x)$ we
set
\[
\|F\|_{W(L^{q_1},L^{q_2})_tW(\Fur
L^{r_1},L^{r_2})_x}= \|
\|F(t,\cdot)\|_{W(\Fur
L^{r_1},L^{r_2})_x}\|_{W(L^{q_1},L^{q_2})_t}.
\]
Observe that (\cite{cordero})
\[
\|F\|_{W(L^{q_1},L^{q_2})_tW(\Fur
L^{r_1},L^{r_2})_x}=
\|F\|_{W\left(L^{q_1}_t(W(\Fur
L^{r_1}_x,L^{r_2}_x)),L^{q_2}_t\right)}.
\]
The following properties of
Wiener amalgam spaces  will
be frequently used in the
sequel.
\begin{lemma}\label{WA}
  Let $B_i$, $C_i$, $i=1,2,3$, be Banach spaces  such that $W(B_i,C_i)$ are well
  defined. Then,
  \begin{itemize}
  \item[(i)] \emph{Convolution.}
  If $B_1\ast B_2\hookrightarrow B_3$ and $C_1\ast
  C_2\hookrightarrow C_3$, we have
  \begin{equation}\label{conv0}
  W(B_1,C_1)\ast W(B_2,C_2)\hookrightarrow W(B_3,C_3).
  \end{equation}
   In particular, for every
$1\leq p, q\leq\infty$, we
have
\begin{equation}\label{p2}
\|f\ast
u\|_{\fpq}\leq\|f\|_{W(\Fur
L^\infty,L^1)}\|u\|_{\fpq}.
\end{equation}
  \item[(ii)]\emph{Inclusions.} If $B_1 \hookrightarrow B_2$ and $C_1 \hookrightarrow C_2$,
   \begin{equation*}
   W(B_1,C_1)\hookrightarrow W(B_2,C_2).
  \end{equation*}
  \noindent Moreover, the inclusion of $B_1$ into $B_2$ need only hold ``locally'' and the inclusion of $C_1 $ into $C_2$  ``globally''.
   In particular, for $1\leq p_i,q_i\leq\infty$, $i=1,2$, we have
  \begin{equation}\label{lp}
  p_1\geq p_2\,\mbox{and}\,\, q_1\leq q_2\,\Longrightarrow W(L^{p_1},L^{q_1})\hookrightarrow
  W(L^{p_2},L^{q_2}).
  \end{equation}
  \item[(iii)]\emph{Complex interpolation.} For $0<\theta<1$, we
  have
\[
  [W(B_1,C_1),W(B_2,C_2)]_{[\theta]}=W\left([B_1,B_2]_{[\theta]},[C_1,C_2]_{[\theta]}\right),
  \]
if $C_1$ or $C_2$ has
absolutely continuous norm.
    \item[(iv)] \emph{Duality.}
    If $B',C'$ are the topological dual spaces of the Banach spaces $B,C$ respectively, and
    the space of test functions $\cC_0^\infty$ is dense in both $B$ and $C$, then
\begin{equation}\label{duality}
W(B,C)'=W(B',C').
\end{equation}
  \end{itemize}
  \end{lemma}
The proof of all these
results can be found in
  (\cite{feichtinger80,feichtinger83,feichtinger90,Heil03}).\par

\section{Fixed-time estimates}\label{section3}
In this section we study estimates for the solution $u(t,x)$ of the Cauchy problem \eqref{cp}, for fixed $t$. We take advantage of the explicit formula for the solution
\begin{equation}\label{sol}
u(t,x)=(K_t\ast u_0)(x),
\end{equation}
where
\begin{equation}\label{chirp0}
K_t(x)=\frac{1}{(4\pi i t)^{d/2}}e^{i|x|^2/(4t)}.
\end{equation}
We already know that
\eqref{chirp0} is in the
Wiener amalgam space $W(\cF
L^1, L^\infty)$ see
\cite{benyi,cordero,baoxiang}.
This is the finest Wiener
amalgam space-norm for
\eqref{chirp0} which,
consequently, gives the worst
behaviour in the time
variable. We aim at improving
the latter, at the expenses
of a rougher $x$-norm.  This
is achieved in the following
result. Indeed, $W(\cF L^1,
L^\infty)\subset W(\cF L^p,
L^\infty) $ for $1\leq p\leq
\infty$ and the case $p=1$
recaptures \cite[Proposition
3.2]{cordero}.

\begin{proposition}\label{chirp}
For $a\in\mathbb{R}$,
$a\not=0$, let $f_{ai}(x)=(a
i)^{-d/2}e^{-\pi
|x|^2/(ai)}$. Then, for
$1\leq p\leq\infty$,
$f_{ai}\in W(\cF L^p,
L^\infty)$ and
\begin{equation}\label{chirpnormp}
\|f_{ai}\|_{W(\cF L^p,
L^\infty)}\asymp
|a|^{-d/p}(1+a^2)^{({d}/{2})(1/p-1/2)}.
\end{equation}
\end{proposition}
\begin{proof} It follows the footsteps
of \cite[Proposition
3.2]{cordero}. We
consider the Gaussian
$g(t)=e^{-\pi |t|^2}$ as
window function to compute $$
\|f_{ai}\|_{W(\cF L^p,
L^\infty)}\asymp
\sup_{x\in\rd}\|\hat{f_{ai}}\ast
M_{-x}g \|_{L^p}.$$ Using
$\widehat{f_{ai}}(\o)=e^{-\pi
a i\o^2}$, we have, for
$p<\infty$,
\begin{align*}
\|\hat{f_a}\ast M_{-x}g \|_{L^p}&=\left(\intrd\left|\intrd e^{-\pi a i(\o-y)^2}e^{-2\pi i x y}e^{-\pi^2 |y|^2}\,dy\right|^p d\o\right)^{1/p} \\
&=(1+a^2)^{-d/4}|a|^{-d/p}\left(\intrd  e^{-\pi p|z|^2/(1+a ^2 )} dz\right)^{1/p}\\
  &=(1+a^2)^{(d/2)(1/p-1/2)}|a|^{-d/p}p^{-d/(2p)}.
  \end{align*}
  Since the right-hand side does not depend on $x$, taking the supremum on $\rd$ with respect to the $x$-variable we attain the desired
  estimate.\par
If $ p =\infty$,
$$\|f_a\|_{W(\cF L^\infty,
L^\infty)} \asymp\sup_{\o\in\rd}|(1+
ai)^{-d/2}e^{-\pi(x-a\o)^2/(1+ai)}|=(1+ a^2)^{-d/4}$$ and we are
done.
\end{proof}

For $a=4\pi t$, $t\not=0$, we infer:
\begin{corollary} Let $K_t(x)$ be the kernel
 in \eqref{chirp0}.
Then, if $1\leq p
\leq\infty$,
\begin{equation}\label{kernelnormp}
\|K_t\|_{W(\cF L^p,
L^\infty)}\asymp
|t|^{-d/p}(1+
t^2)^{(d/2)(1/p-1/2)}.\end{equation}
\end{corollary}
\begin{lemma}
 Let $1\leq p,q,r\leq \infty$, with
$1/r=1/p+1/q$, then
\begin{equation}\label{eq01}W(\Fur L^p,L^\infty)\ast W(\Fur
L^q,L^1)\hookrightarrow W(\Fur L^r,L^\infty).
\end{equation}
\end{lemma}
\begin{proof}
This is a consequence of the convolution relations for Wiener
amalgam spaces in Lemma \ref{WA} (i), being $\Fur L^p\ast \Fur
L^q= \Fur (L^p\cdot L^q)=\Fur L^r$ by
 H\"older's Inequality with
$1/r=1/p+1/q$, and $L^\infty\ast L^1 \hookrightarrow L^\infty$.
\end{proof}

\begin{prop}
It turns out, for $2\leq
q\leq\infty$,
\begin{equation}\label{est1}
 \|e^{it\Delta}u_0\|_{W(\Fur L^{q^\prime},
 L^\infty)}\lesssim
|t|^{d(2/q-1)}(1+t^2)^{d(1/4-1/q)}
\|u_0\|_{W(\Fur L^q,L^1)}.
\end{equation}
\end{prop}
\begin{proof} We use the explicit
representation of the Schr\"odinger evolution operator $ e^{it\Delta}u_0(x)=(K_t\ast
u_0)(x)$.  Let $1\leq p, q\leq\infty$, satisfying
\begin{equation}\label{har}\frac1p+\frac2q=1.
\end{equation}
Then the kernel norm
\eqref{kernelnormp} and  the
convolution relations
\eqref{eq01} yield the
desired result.
\end{proof}
\begin{theorem}\label{tfix}
 For $2\leq q, r,s\leq
\infty$ such that
\begin{equation}\label{indrel}
\frac1 s=\frac1 r+\frac 2
q\left(\frac12-\frac1r\right),
\end{equation}
we have
\begin{equation}
\|e^{it\Delta}u_0\|_{W(\Fur
L^{s^\prime},
 L^r)}\lesssim |t|^{d\left(\frac2q-1\right)
\left(1-\frac2r\right)}(1+t^2)^{d\left
(\frac{1}{4}
 -\frac{1}{q}\right)\left(1-\frac{2}{r}\right)}\nonumber\\
\|u_0\|_{W(\Fur
L^{s},L^{r^\prime})}.\label{est2n}
\end{equation}
\end{theorem}
In particular, for $2\leq
r\leq \infty$,
\begin{equation}\label{est2i}
\|e^{it\Delta}u_0\|_{W(
L^{2}, L^r)}\lesssim (1+
t^2)^{-\frac{d}2\left(\frac{1}{2}
-\frac{1}{r}\right)}\|u_0\|_{W(L^{2},
L^{r^\prime})}.
\end{equation}
\begin{proof}
Estimates  \eqref{est2n}
follow by complex
interpolation of
\eqref{est1}, which
corresponds to $r=\infty$,
with the $L^2$ conservation
law
\begin{equation}\label{l2}
\|e^{it\Delta}u_0\|_{L^2}=\|u_0\|_{L^2},
\end{equation}
which corresponds to $r=2$.\\
Indeed, $L^2=W(\Fur
L^2,L^2)=W( L^2,L^2)$. By
Lemma \ref{WA}, item (iii),
for $0<\theta=2/r<1$, and
$1/{s^\prime}=(1-2/r)/q^{\prime}+1/r$,
so that relation
\eqref{indrel} holds,
\[
\left[W(\Fur L^{q^\prime}, L^\infty),W(\Fur L^2,L^2)\right]_{[\theta]}=W\left([\Fur L^{q^\prime},\Fur L^2]_{[\theta]}, [L^\infty,L^2]_{[\theta]}\right)=W(\Fur L^{s^\prime}, L^{r})
\]
and
\[
\left[W(\Fur L^q, L^1),W(\Fur L^2,L^2)\right]_{[\theta]}=W\left([\Fur L^q,\Fur L^2]_{[\theta]}, [L^1,L^2]_{[\theta]}\right)=W(\Fur L^{s}, L^{r^\prime}),
\]
so that the estimate
\eqref{est2n} is attained.
\end{proof}

\section{Strichartz estimates}\label{section4}
In this section we prove
Theorem \ref{prima2i}.
Precisely, we first study the
case
$q_1=\tilde{q}_1=\infty$,
$r_1=\tilde{r}_1=2$. In view
of the inclusion relation of
Wiener amalgam spaces, it
suffices to prove it for the
pairs $(q_2,r_2)$ scale
invariant (i.e. satisfying
\eqref{pri2} with equality),
namely the following result.
\begin{theorem}\label{prima2}
Let $2\leq q\leq\infty$,
$2\leq r\leq\infty$, such
that
\begin{equation*}\frac{2}{q}+\frac{d}{r}
=\frac{d}{2},
\end{equation*}
$(q,r,d)\not=(2,\infty,2)$,
and similarly for
$\tilde{q},\tilde{r}$. Then
we have the homogeneous
Strichartz estimates
\begin{equation}\label{hom2}\|e^{it\Delta} u_0\|_{W(L^{\infty},L^q)_t W(
L^{2},L^r)_x}\lesssim \|u_0\|_{L^2_x},
\end{equation}
the dual homogeneous Strichartz estimates
\begin{equation}\label{dh2}
\|\int e^{-is\Delta} F(s)\,ds\|_{L^2}\lesssim
\|F\|_{W(L^{1},L^{\tilde{q}^\prime})_t W(
L^{2},L^{\tilde{r}^\prime})_x},
\end{equation}
and the retarded Strichartz estimates
\begin{equation}\label{ret2}
\|\int_{s<t} e^{i(t-s)\Delta} F(s)\,ds\|_{W(L^\infty,L^q)_t W(
L^{2},L^r)_x}
\lesssim\|F\|_{W(L^{1},L^{\tilde{q}^\prime})_t
W(L^{2},L^{\tilde{r}^\prime})_x}.
\end{equation}
\end{theorem}
\begin{proof}[Proof in the non-endpoint case]
 Here we prove Theorem \ref{prima2} in the
  non-endpoint case, namely for
  $q>2,\tilde{q}>2$.
 The techniques are quite similar to those
  in \cite[Theorem 1.1]{cordero}.
We shall sketch the proof of
the homogeneous and dual
Strichartz estimates. We
first  show the estimate
\eqref{hom2}. The case:
$q=\infty$, $r=2$, follows at
once from the conservation
law \eqref{l2}. Indeed,
$W(L^{\infty},L^\infty)_t=L^\infty_t$
and $ W( L^{2},L^2)_x=L^2_x$,
so that,  taking the supremum
over $t$ in $\|e^{it\Delta}
u_0\|_{L^2_x}=\|u_0\|_{L^2_x}$,
we attain the claim.

 To prove the remaining cases, we can apply the usual
$TT^\ast$ method (or ``orthogonality principle", see \cite[Lemma 2.1]{GinibreVelo92} or \cite[page 353]{stein}),
 because of H\"older's type inequality
\begin{equation}\label{holder}
|\langle F,G\rangle_{L^2_t L^2_x}|\leq \|F\|_{W(L^\infty,L^q)_t W(
L^{2},L^r))_x}\|G\|_{W(L^{1},L^{q^\prime})_t W(
L^{2},L^{r^\prime})_x},
\end{equation}
which can be proved directly from the definition of these spaces.\\
As a consequence, it suffices to prove the estimate
\begin{equation}\label{aus}
\|\int e^{i(t-s)\Delta}F(s)\,ds\|_{W({L^{\infty}},{L^q})_t W(
L^{2},L^r)_x}\lesssim\|F\|_{W(L^{1},L^{q^\prime})_t
W(L^2,L^{r^\prime})_x}.
\end{equation}

Recall the
Hardy-Littlewood-Sobolev
fractional integration
theorem (see e.g.
\cite{stein}, page 119) in
dimension $1$:
\begin{equation}\label{conv1}
L^p(\bR)\ast L^{1/\alpha,\infty}(\bR)\hookrightarrow L^q(\bR),
\end{equation}
for $1\leq p<q<\infty$,
$0<\alpha<1$, with
$1/p=1/q+1-\alpha$ (here
$L^{1/\alpha,\infty}$ is the
weak $L^{1/\alpha}$ space,
see e.g. \cite{steinweiss}).
 Now, set
$\alpha=d(1/2-1/r)=2/q$
($q>2$) so that, for
$p=q^\prime$,
$L^{q^\prime}\ast
L^{1/\alpha,\infty}
\hookrightarrow L^q$.
Moreover, observe that
$(1+|t|)^{-\a}\in
W(L^\infty,L^{1/\a,\infty})(\bR)$.
The convolution relations
\eqref{conv0} then give
$$W(L^1,L^{q^\prime})(\bR) \ast W(L^\infty,L^{1/\a,\infty})(\bR) \hookrightarrow W(L^\infty,L^{q})(\bR).
$$
The preceding relations, together with the fixed-time estimates
\eqref{est2i} and Minkowski's inequality allow to write
\begin{align*}
\|\int e^{i(t-s)\Delta}F(s)\,ds&\|_{W(L^{\infty},{L^q})_t W(
L^{2},L^r)_x}\\
&\leq \|\int\|e^{i(t-s)\Delta} F(s)\|_{W(
L^{2},L^r)_x}\,ds\|_{W(L^{\infty},{L^q})_t}\\
&\lesssim\|\int\|F(s)\|_{W(L^{2},L^{r^\prime})_x}
(1+|t-s|)^{-\a} ds\|_{W(L^\infty,{L^q})_t}\\
&\lesssim \|F\|_{W(L^{1},L^{q^\prime})_t
W(L^2,L^{r^\prime})_x}.
\end{align*}
The estimate \eqref{dh2} follows  from \eqref{hom2} by duality.
\par
The same arguments as in \cite[page 13]{cordero} then give the
retarded Strichartz estimate \eqref{ret2}.
\end{proof}

\begin{proof}[Proof in the
endpoint
case]\label{section5}

We are going to prove
Theorem \ref{prima2} in the
endpoint case
$(q,r)=P:=(2,2d/(d-2))$ or
$(\tilde{q},\tilde{r})=P$,
$d>2$. Hence, we prove the
estimates
\begin{equation}\label{hombis}\|e^{it\Delta}
 u_0\|_{W(L^{\infty},L^2)_t W(
L^{2},L^r)_x}\lesssim
\|u_0\|_{L^2_x},\quad
r=2d/(d-2),
\end{equation}
\begin{equation}\label{dhbis}
\|\int e^{-is\Delta}
F(s)\,ds\|_{L^2}\lesssim
\|F\|_{W(L^{1},L^{\tilde{q}^\prime})_t
W(
L^{2},L^{\tilde{r}^\prime})_x},\quad
\tilde{r}=2d/(d-2),
\end{equation}
and
\begin{equation}\label{ret2bis}
\|\int_{s<t} e^{i(t-s)\Delta}
F(s)\,ds\|_{W(L^\infty,L^2)_t
W( L^{2},L^r)_x}
\lesssim\|F\|_{W(L^{1},L^{2})_t
W(L^{2},L^{\tilde{r}^\prime})_x},
\end{equation}
with $(q,r)$,
$(\tilde{q},\tilde{r})$
Schr\"odinger admissible, and
$(q,r)=P$ or $(\tilde
{q},\tilde{r})=P$.\par We
follow the pattern in
\cite{cordero,keel}, with
several changes according to
our setting. Hence, we study
bilinear form estimates via a
Withney decomposition in time
(see \eqref{dec}). We
estimate each dyadic
contribution in
\eqref{primolemma} and
\eqref{lammasecondo}. Finally
we conclude by a lemma of
real interpolation theory to
sum these estimates.
\par
Precisely, by the same
duality arguments as the ones
used in the previous part of
the proof, we observe that it
suffices to prove
\eqref{hombis}. This is
equivalent to the bilinear
estimate
\[
|\iint \langle e^{-is\Delta} F(s),
e^{-it\Delta}G(t)\rangle\,ds\,dt|\lesssim \|F\|_{W(L^1,L^{2})_t W(
L^{2},L^{{r}^\prime})_x} \|G\|_{W(L^1,L^{2})_t
W(L^2,L^{{r}^\prime})_x}.
\]
By symmetry, it is enough to prove
\begin{equation}\label{aim}
|T(F,G)|\lesssim \|F\|_{W(L^1,L^{2})_t W( L^{2},L^{{r}^\prime})_x}
\|G\|_{W(L^1,L^{2})_t W( L^{2},L^{{r}^\prime})_x},
\end{equation}
where
\[
T(F,G)=\iint_{s<t} \langle e^{-is\Delta} F(s),
e^{-it\Delta}G(t)\rangle\,ds\,dt.
\]
Here the critical exponent
$q=2$ appears in the global
component, which control the
 \emph{decay in the $t$-variable at
infinity}, hence the form
$T(F,G)$ is decomposed
dyadically as
\begin{equation}\label{dec} T=\tilde{T}+\sum_{j\geq 0} T_j,
\end{equation}
with
\begin{equation}\label{ttilde}
\tilde{T}(F,G)=\iint_{t-1<s< t} \langle e^{-is\Delta} F(s),
e^{-it\Delta}G(t)\rangle\,ds\,dt
\end{equation}
and
\begin{equation}\label{tj}
T_j(F,G)=\iint_{t-2^{j+1}<s\leq t-2^j} \langle e^{-is\Delta} F(s),
e^{-it\Delta}G(t)\rangle\,ds\,dt.
\end{equation}
In the sequel we shall study the behaviour of $\tilde{T}$ and $T_j$
separately. We shall use repeatedly the following fact.
\begin{lemma}\label{fcompact} Let
$1\leq  p,q,r\leq\infty$ be such that
 $1/p+1/q=1/r$. For every $f\in W(L^1,L^q)(\bR)$
supported in any interval $I$
 of length
$L\geq1$, it turns out
\begin{equation}\label{normcomp}
    \|f\|_{W(L^1,L^r)}\leq C_p L^{1/p}\|f\|_{W(L^1,L^q)}.
\end{equation}
\end{lemma}
\begin{proof}
To compute the
$W(L^1,L^r)$ norm of $f$ we
choose $g=\chi_{[0,1]}$, the
characteristic function of
the interval $[0,1]$. Then,
$\|f\|_{W(L^1,L^r)}\asymp
\|\,\|
fT_yg\|_{L^1}\|_{L^r}$. Since
$f$ is supported in an
interval $I$ of length $L$,
  the mapping
$$y\longmapsto\|
fT_yg\|_{L^1}
$$
is supported in an interval $\tilde{I}$ of length $L+2$. Hence, for
$1/p+1/q=1/r$, H\"older's inequality yields
\begin{eqnarray*}
\| \| fT_yg\|_{L^1}\|_{L^r}
&=& \left(\int_{\tilde{I}}\|
fT_yg\|_{L^1}^r dy
\right)^{1/r}\leq
(L+2)^{1/p}\|\,\|
fT_yg\|_{L^1}\|_{L^q} \\
&\leq & C_p L^{1/p}\|f\|_{W(L^1,L^q)}
\end{eqnarray*}as desired.
\end{proof}
\begin{lemma}
We have
\begin{equation}\label{primolemma}
|\tilde{T}(F,G)|\lesssim \|F\|_{W(L^1,L^{2})_t W(
L^{2},L^{{r}^\prime})_x} \|G\|_{W(L^1,L^{2})_t W(
L^{2},L^{{r}^\prime})_x}.
\end{equation}
\end{lemma}
\begin{proof}
We assume $F$ and $G$
compactly supported, with
respect to the time variable,
in intervals of duration $1$
(indeed, in \eqref{ttilde},
$F$ and $G$ can be replaced
by $\chi_{[t-1,t]}(s) F(s)$
and $\chi_{[s,s+1]}(t) G(t)$,
respectively).

Since $|t-s|\leq1$, it
follows from the duality
properties \eqref{duality}
and the fixed-time estimate
\eqref{est2i} that
\begin{eqnarray*}
|\langle e^{-is\Delta} F(s), e^{-it\Delta}G(t)\rangle|&=&
|\langle  F(s), e^{i(s-t)\Delta}G(t)\rangle|\\
&\lesssim&\|F(s)\|_{W( L^{2},L^{{r}^\prime})_x}
\|e^{i(s-t)\Delta}G(t) \|_{W(
L^{2},L^{{r}})_x}\\
&\lesssim& \|F(s)\|_{W( L^{2},L^{{r}^\prime})_x}(
1+|s-t|)^{-d\left(\frac12-\frac1r\right)} \|G(t)\|_{W(L^2
,L^{{r}^\prime})_x}\\
&\lesssim& \|F(s)\|_{W( L^{2},L^{{r}^\prime})_x} \|G(t)\|_{W(L^2
,L^{{r}^\prime})_x}
\end{eqnarray*}
Integrating with respect to the variables $s$ and $t$ we then obtain
$$|\tilde{T}(F,G)| \lesssim \|F\|_{L^1_tW( L^{2},L^{{r}^\prime})_x} \|G\|_{L^1_tW(L^2
,L^{{r}^\prime})_x}.
$$
Lemma \ref{fcompact} with
$p=q=2$ and $r=1$, applied to
each function $\|F(t)\|_{W(
L^{2},L^{{r}^\prime})_x}$ and
$\|G(t)\|_{W(
L^{2},L^{{r}^\prime})_x}$,
gives the result.
\end{proof}\\
The estimates of the pieces
$T_j(F,G)$ follow the
techniques of \cite[Lemma
4.1]{keel} and \cite[Lemma
5.2]{cordero}, adapted to our
context.
\begin{lemma} We have
\begin{equation}\label{lammasecondo}
|T_j(F,G)|\lesssim 2^{-j\beta(a,b)}\|F\|_{W(L^1,L^2)_t W(
L^2,L^{a^\prime})}\|G\|_{W(L^1,L^2)_t W(L^2,L^{b^\prime})},
\end{equation}
for $(1/a,1/b)$ in a neighborhood of $(1/r,1/r)$.
\end{lemma}
\begin{proof}
Observe that here $d>2$ hence $r<\infty$. Then, the  result follows
by complex interpolation (Lemma \ref{WA}, (iii)) from the following
cases:
\begin{itemize}
\item[(i)] $a=\infty$, $b=\infty$, \item[(ii)] $2\leq a<r$, $b=2$,
\item[(iii)] $a=2$, $2\leq b<r$.
\end{itemize}
{\it Case (i)}. We need to show the estimate
\begin{equation}\label{(i)}|T_j(F,G)|\lesssim 2^{-j(\frac {d}2-1)} \|F\|_{W(L^1,L^2)_t W(
L^2,L^1)}\|G\|_{W(L^1,L^2)_t W(L^2,L^1)}.
\end{equation}
Here the fixed-time estimate \eqref{est2i} and $|t-s|\lesssim 2^{j}$
yield
\[
|\langle e^{-is\Delta} F(s), e^{-it\Delta} G(t)\rangle|\lesssim
2^{-j\frac{d}2}\|F(s)\|_{W(L^2,L^1)}\|G(t)\|_{W(L^2,L^1)}.
\]
Integrating with respect to the variables $s$ and $t$,
\[
|T_j(F,G)|\lesssim
2^{-j\frac{d}2}\|F\|_{L^1_t(W(L^2,L^1)_x)}
\|G\|_{L^1_t(W(L^2,L^1)_x)}.
\]
Again, we can assume $F$ and
$G$ compactly supported, with
respect to the time variable,
in intervals of duration $
2^j$. Applying Lemma
\ref{fcompact} with $p=q=2$
and $r=1$,  to both functions
$\|F(t)\|_{W(L^2,L^1)_x}$ and
$\|G(t)\|_{W(L^2,L^1)_x}$ we
attain
\eqref{(i)}.\noindent\par
 {\it Case (ii)}. We have to show
\begin{equation}\label{(ii)}|T_j(F,G)|\lesssim 2^{-j(\frac d4-1-\frac d {2a})} \|F\|_{W(L^1,L^2)_t W(L^2,L^{a^\prime})_x}
\|G\|_{W(L^1,L^2)_t  L^2_x}.
\end{equation}
Using similar arguments to the previous case we obtain
\begin{equation}\label{unozero}
|T_j(F,G)|\lesssim \sup_{t}\|\int_{t-2^{j+1}<s\leq t-2^j}
e^{-is\Delta}F(s)\,ds\|_{L^2_x}\|G\|_{L^1_tL^2_x},
\end{equation}
and \begin{equation}\label{unouno}\|G\|_{L^1_tL^2_x}\lesssim
2^{j/2}\|G\|_{W(L^1,L^2)_tL^2_x}.
\end{equation}
For $a\geq 2$, let now $\tilde{q}=\tilde{q}(a)$ be defined by
\begin{equation}\label{qa}
\frac{2}{\tilde{q}(a)}+\frac{d}{a}=\frac{d}{2}.
\end{equation}
The non-endpoint case of \eqref{dh2}, written for $\tilde{r}=a$ and
the $\tilde{q}$ above, gives
\begin{eqnarray*}\sup_t\!\|\int_{t-2^{j+1}<s\leq t-2^j}\!\!\!
e^{-is\Delta}F(s)\,ds\|_{L^2_x}\!\!&=&\!\!\sup_t\|\int_{\R}
e^{-is\Delta}(T_{-t}\chi_{[-2^{j+1},-2^{j}]})(s)F(s)\,ds\|_{L^2_x}\\
&\lesssim& \|F\|_{W(L^1,L^{\tilde{q}^\prime})W(L^2,L^{a'})},
\end{eqnarray*}
for every $2\leq a<r$. Since the support of $F$ with respect to the
time is contained in an interval of
 duration $2^j$, Lemma
\ref{fcompact} with $q=2$, $r=\tilde{q}^\prime$, and
$$\frac{1}{p}=\frac{1}{\tilde{q}^\prime}-\frac{1}{2}=\frac{1}{2}-\frac{1}{\tilde{q}}=\frac12-\frac
d4+\frac d{2a},
$$
gives
$$\|F\|_{W(L^1,L^{\tilde{q}^\prime})W(L^2,L^{a'})}\lesssim 2^{j(\frac12+\frac d{2a}-\frac d4)}
\|F\|_{W(L^1,L^2)W(L^2,L^{a'})}.
$$
This estimate, together with
\eqref{unozero} and
\eqref{unouno}, yields the
estimate
\eqref{(ii)}.\par\noindent
{\it Case (iii).} Use the
same arguments as in case
({\it ii}).
\end{proof}

It remains to show
\begin{equation}\label{abbnm}
\sum_{j\geq0}|T_j(F,G)|\lesssim \|F\|_{W(L^1,L^2)_t
W(L^2,L^{r^\prime})_x} \|G\|_{W(L^1,L^2)_t W(L^2,L^{{r}^\prime})_x}.
\end{equation}
Now, \eqref{abbnm} can be achieved from
\eqref{lammasecondo} and some real
interpolation results (collected in
Appendix A below), as in
\cite{cordero,keel}.\par In details, we
single out $a_0,a_1,b_0,b_1$ such that
$(1/r,1/r)$ is inside a small triangle
with vertices $(1/a_0,1/b_0)$,
$(1/a_1,1/b_0)$ and $(1/a_0,1/b_1)$
(see Figure 2), so that
\[
\beta(a_0,b_1)=\beta(a_1,b_0)\not=\beta(a_0,b_0).
\]
\vspace{-0.2cm}
  \begin{center}
           \includegraphics{figstr2.1}
            \\
           $ $
\end{center}
           \begin{center}{Figure 2.  }
           \end{center}
  \vspace{1.2cm}
Then, we apply Lemma
\ref{interp} with\footnote{
We set
$l^s_q=L^q(\mathbb{Z},2^{js}dj)$,
where $dj$ is the counting
measure. Recall (see
\cite[Section
5.6]{bergh-lofstrom}) that
$(l_\infty^{s_0},l^{s_1}_\infty)_{\theta,1}
=l^s_1$ whenever
$s_0\not=s_1$ and
$s=(1-\theta)s_0+\theta s_1$,
$0<\theta<1$.} $T=\{T_j\}$,
(after setting $T_j=0$ for
$j<0$)
$C_0=l_\infty^{\beta(a_0,b_0)}$,
$C_1=l_\infty^{\beta(a_0,b_1)}$
and,  for $k=0,1$, we take
$$A_k=W(L^1(W(L^2,L^{{a_k}^\prime})_x),L^{2}),
\  B_k=W(L^1(W(L^2,L^{{b_k}^\prime})_x),L^{2}). $$
 Here we choose
$\theta_0$, $\theta_1$, so that
$$1/r=(1-\theta_0)/a_0+\theta_0/a_1, \quad
1/r=(1-\theta_1)/b_0+\theta_1/b_1.$$
The assumptions are satisfied
in view of
\eqref{lammasecondo}.
Moreover, with
$\theta=\theta_0+\theta_1$,
we have
$$
(1-\theta)\beta(a_0,b_0)+\theta\beta(a_0,b_1)
=\beta(r,r)=0. $$
 Hence we
attain the desired estimate
\eqref{abbnm}, because
\begin{eqnarray*}
  (A_0,A_1)_{\theta_0,2} &\hookleftarrow& W((L^1(W(L^2,L^{{a_0}^\prime})_x),L^1(W(L^2,L^{{a_1}^\prime})_x))_{\theta_0,2},L^{2})\\
   &\hookleftarrow& W(L^1( W(L^2,L^{{a_0}^\prime})_x,W(L^2,L^{{a_1}^\prime})_x)_{\theta_0,2}),L^2)\\
   &\hookleftarrow& W(L^1( W(L^2,L^{{r}^\prime})_x),L^2)=W(L^1,L^2)_t
   W(L^2,L^{{r}^\prime})_x,
\end{eqnarray*}
where we used Proposition
\ref{inter9} for the first
and third embedding and
Proposition \ref{interlp} for
the second one (the same
holds for
$(B_0,B_1)_{\theta_0,2}$).
\par Similarly to \cite[page
19]{cordero} one can obtain
the retarded estimates. This
concludes the proof of
Theorem \ref{prima2}.
\end{proof}\par\medskip\noindent
The Strichartz estimates in Theorem \ref{prima2}
 can be
combined with the classical ones
\eqref{S1} to  obtain the estimates in
Theorem \ref{prima2i}.\par
\medskip\noindent
{\textbf  {Proof of Theorem \ref{prima2i}.}} We prove the homogeneous estimates \eqref{hom2i}, the other ones follow by similar arguments.
 Recall that
the classical homogeneous
Strichartz estimates can be
written as
\begin{equation}\label{S1Wbis}
\|e^{it\Delta}
u_0\|_{W(L^{\tilde{q}},L^{\tilde{q}})_t
W(L^{\tilde{r}},L^{\tilde{r}})_x}\lesssim
\|u_0\|_{L^2_x},
\end{equation}
for every Schr\"odinger
admissible pair
$(\tilde{q},\tilde{r})$, i.e.
$\tilde{q},\tilde{r}\geq2$,
with
$2/\tilde{q}+d/\tilde{r}=d/2$,
$(\tilde{q},\tilde{r},d)\not=(2,\infty,2)$.
By complex interpolation
between \eqref{S1Wbis} and
\eqref{hom2}
 one has
\begin{equation}\label{S1New}
\|e^{it\Delta} u_0\|_{W(L^{q_1},L^{q_2})_t W(
L^{r_1},L^{r_2})_x}\lesssim \|u_0\|_{L^2_x},
\end{equation}
with  $(q_i,r_i)$, $i=1,2$,
Schr\"odinger admissible
pairs.\\
 Here
$$\frac1{q_1}=\frac{1-\theta}{\tilde{q}}+\frac\theta \infty,
\quad \frac1{q_2}=\frac{1-\theta}{\tilde{q}}+\frac\theta
 {q},
$$
so that
$$\frac1{q_2}=\frac1{q_1}+\frac\theta
{q}.
$$
Hence $q_1\geq q_2$, i.e.,
$r_1\leq r_2$. This shows
\eqref{hom2} when $(q_1,r_1)$
and $(q_2,r_2)$ satisfy
\eqref{pri1} and \eqref{pri2}
with equality. The general
case follows from the
inclusion relations of Wiener
amalgam spaces, which allow
us to increase $q_2,r_2$, and
diminish $q_1,r_1$ (see
Figure 1).\par

\section{Sharpness of
fixed-time and Strichartz
estimates}\label{necessarie}
In this section we prove the
sharpness of the estimates
\eqref{est2}, \eqref{str} and
\eqref{hom2i}. To this end we
need the following three
lemmata.\par We recall from
\cite[page 257]{folland} the
following well-known formula
for Gaussian integrals.
\begin{lemma}
Let $A$ be a $d\times d$
complex matrix such that
$A=A^\ast$ and ${\rm Re}\,A$
is positive definite. Then
for every $z\in\C^n$,
\begin{equation}\label{lem0}
\int e^{-\pi x A x-2\pi i z
x}dx=({\rm det}\,{A})^{-1/2}
e^{-\pi zA^{-1}z},
\end{equation}
where the branch of the
square root is determined by
the requirement that $({\rm
det}\,{A})^{-1/2}>0$ when $A$
is real and positive
definite.
\end{lemma}
\begin{lemma}\label{lem1}
For $c\in\C$, $c\not=0$,
consider the function
$\phi^{(c)}(x)=e^{-\pi
c|x|^2}$, $x\in\R^d$. For
every $c_1,c_2\in\C$, with ${\rm
Re}\, c_1\geq0$, ${\rm Re}\,
c_2>0$, we have
\begin{equation}
\phi^{(c_1)}\ast\phi^{(c_2)}=(c_1+c_2)^{-d/2}
\phi^{(\frac{c_1c_2}{c_1+c_2})}.
\end{equation}
\end{lemma}
\begin{proof}
Using the equality \eqref{lem0},
\begin{align}
\left(\phi^{(c_1)}\ast\phi^{(c_2)}\right)(x)
&=\int e^{-\pi c_1|x-t|^2-\pi
c_2|t|^2}dt\nonumber\\&=e^{-\pi
c_1|x|^2}\int
e^{-\pi(c_1+c_2)|t|^2-2\pi
c_1 x t}dt\nonumber\\
&=(c_1+c_2)^{-d/2} e^{-\pi
c_1|x|^2}e^{\frac{\pi
c_1^2|x|^2}{c_1+c_2}},\nonumber\\
&=(c_1+c_2)^{-d/2}
\phi^{(\frac{c_1c_2}{c_1+c_2})}(x),\nonumber
\end{align}
as desired.
\end{proof}\\
In particular,
the solution of the Cauchy
problem \eqref{cp}, with
initial datum
$u_0(x)=e^{-\pi|x|^2}$, is
given by the formula
\begin{equation}\label{zz1}
u(t,x)=(1+4\pi i t)^{-d/2}
e^{-\frac{\pi|x|^2}{1+4\pi i
t}}.
\end{equation}
This follows at once from
Lemma \ref{lem1}, since the solution $u(t,x)$ can be rephrased as
\[
u(t,x)=(K_t\ast
u_0)(x)=\frac{1}{(4\pi i
t)^{d/2}} \left(\phi^{((4\pi
i t)^{-1})}\ast
\phi^{(1)}\right)(x).
\]

\begin{lemma}\label{lemm2}
For $a,b\in\R$, $a>0$, set
$f_{a+ib}(x)=(a+ib)^{-d/2}e^{-\frac{\pi|x|^2}{a+ib}}$.
Then, for every $1\leq
q,r\leq\infty$,
\begin{equation}\label{lem2}
\|f_{a+ib}\|_{W(\Fur
L^q,L^r)}\asymp \frac{\left((a+1)^2+b^2\right)^{
\frac{d}{2}\left(\frac{1}{q}-\frac{1}{2}
\right)}}{a^\frac{d}{2r}\left(a(a+1)+b^2\right)^{\frac{d}{2}\left(\frac{1}{q}-\frac{1}{r}\right)}}.
\end{equation}
\end{lemma}
\begin{proof}  We use the Gaussian  $g(y)=e^{-\pi |y|^2}$ as window function, so that
the Wiener amalgam norm
$W(\Fur L^q,L^r)$ reads
\[
\|f_{a+ib}\|_{W(\Fur
L^q,L^r)}\asymp\|\|f_{a+ib} T_x
g\|_{\Fur L^q}\|_{L^r_x}.
\]
Now,
\begin{align}
\widehat{f_{a+ib} T_x
g}(\omega)&=\left(\widehat{f_{a+ib}}\ast
M_{-x} g\right)(\omega)\\
&=\int
e^{-\pi(a+ib)|\omega-y|^2}
e^{-2\pi ix y}
e^{-\pi|y|^2}dy\nonumber\\
&=e^{-\pi(a+ib)|\omega|^2}\int
e^{-\pi(a+1+ib)|y|^2-2\pi
i(x+i(a+ib)\omega)y}dy\nonumber\\
&=(a+1+ib)^{-d/2}e^{-\pi(a+ib)|\omega|^2}
e^{-\pi\frac{(x+i(a+ib)\omega)
(x+i(a+ib)\omega)}{a+1+ib}},\nonumber
\end{align}
where we used \eqref{lem0}.
Hence, after a simple
computation,
\[
|\widehat{f_{a+ib} T_x
g}(\omega)|=((a+1)^2+b^2)^{-d/4}
e^{-\frac{\pi}{(a+1)^2+b^2}[(a(a+1)+b^2)|
\omega|^2+2bx\omega+(a+1)|x|^2]}.
\]
It follows that
\begin{align}
\|f_{a+ib} T_x g\|_{\Fur
L^q}&=((a+1)^2+b^2)^{-d/4}
e^{-\frac{\pi(a+1)|x|^2}{(a+1)^2+b^2}}\left(\int
 e^{-\frac{q\pi}{(a+1)^2+b^2}[(a^2+b^2+a)|\omega|^2+2bx\omega]}d\omega\right)^{\frac{1}{q}}\nonumber\\
 &=((a+1)^2+b^2)^{-d/4} e^{-\frac{\pi
 a|x|^2}{a(a+1)+b^2}}\left(\int
 e^{-\frac{q\pi}{(a+1)^2+b^2}|
 \sqrt{a(a+1)+b^2}\omega+
 \frac{b}{\sqrt{a(a+1)+b^2}}x|^2}d\omega
 \right)^\frac{1}{q}\nonumber\\
 &=((a+1)^2+b^2)^{-d/4}
e^{-\frac{\pi
a|x|^2}{a(a+1)+b^2}}(a(a+1)+b^2)^{-\frac{d}{2q}}\left(\frac{(a+1)^2+b^2}{q}\right)^{\frac{d}{2q}}\nonumber\\
&=\frac{\left((a+1)^2+b^2\right)^{\frac{d}{2}\left(\frac{1}{q}-\frac{1}{2}\right)}}{q^{\frac{d}{2q}}\left(a(a+1)+b^2\right)^{\frac{d}{2q}}}
e^{-\frac{\pi
a|x|^2}{a(a+1)+b^2}}.\nonumber
 \end{align}
 By taking the $L^r$ norm of
 this expression one obtains
 \eqref{lem2}.
\end{proof}
\begin{proposition}{\rm (Sharpness of
\eqref{est2}).}\label{prop1}
Suppose that, for some fixed
$t_0\in\R$, $1\leq
r\leq\infty$, $C>0$, the
following estimate holds:
\begin{equation}\label{z01}
\|e^{it_0\Delta}
u_0\|_{W(\Fur
L^{r'},L^r)}\leq C
\|u_0\|_{W(\Fur
L^{r},L^{r'})},\quad \forall
u_0\in\mathcal{S}(\R^d).
\end{equation}
Then $r\geq2$.\par
 Assume now
that, for some $\alpha\in\R$,
$C>0$, $\delta>0$, $1\leq
r\leq\infty$, the estimate
\begin{equation}\label{z1}
\|e^{it\Delta} u_0\|_{W(\Fur
L^{r'},L^r)}\leq C
t^\alpha\|u_0\|_{W(\Fur
L^{r},L^{r'})},\quad \forall
u_0\in\mathcal{S}(\R^d),
\end{equation}
holds for every
$t\in(0,\delta)$. Then
\begin{equation}\label{z2}
\alpha\leq
-2d\left(\frac{1}{2}-\frac{1}{r}\right).
\end{equation}
\indent Finally, if
$u_0(x)=e^{-\pi|x|^2}$ then
\begin{equation}\label{z3}
\|e^{it\Delta}u_0\|_{W(\Fur
L^{r'},L^r)}\sim t^{
-d\left(\frac{1}{2}-\frac{1}{r}\right)},\quad
{\it as}\ t\to+\infty.
\end{equation}
\end{proposition}
\begin{proof}
We consider the one parameter
family of initial data
$u_0(\lambda
x)=e^{-\pi\lambda^2|x|^2}$,
$\lambda>0$. If $u$ is the
function in \eqref{zz1}, the
corresponding solutions will
be
\begin{align}\label{c1}
u(\lambda^2 t,\lambda
x)&=(1+4\pi i t
\lambda^2)^{-d/2}
e^{-\frac{\pi\lambda^2|x|^2}{1+4\pi
i t\lambda^2}}\\
&=\lambda^{-d}f_{\lambda^{-2}+4\pi
i t}(x),\nonumber
\end{align}
where we used the notation in
Lemma \ref{lemm2}.\par Now,
\eqref{lem2} gives
\begin{equation}\label{z02}
\|u_0(\lambda\,\cdot)\|_{W(\Fur
L^{r},L^{r'})}=\lambda^{-d}\|f_{\lambda^{-2}}
\|_{W(\Fur L^{r},L^{r'})}\sim
c\lambda^{-d/r'}
\end{equation}
 both for $\lambda\to0^+$ and
 $\lambda\to+\infty$, for some $c>0$.
 On the other hand, again an application of \eqref{lem2}
 yields
 \begin{align}\label{c2}
 \|u(\lambda^2t,\lambda\,\cdot)\|_{
 W(\Fur
 L^{r'},L^{r})}&=\lambda^{-d}\|f_{\lambda^{-2}+4\pi
 i t}\|_{
 W(\Fur
 L^{r'},L^{r})}\nonumber\\
 &\asymp\frac{\lambda^{-d/r'}\left[(1+
 \lambda^{-2})^2+
 t^2\right]^{\frac{d}{2}
 \left(\frac{1}{2}-\frac{1}{r}\right)}}
 {\left[
 \lambda^{-2}(\lambda^{-2}+1)+ t^2
 \right]^{d\left(\frac{1}{2}-
 \frac{1}{r}\right)}}.
 \end{align}
 Now, for fixed $t=t_0$, the
 expression in \eqref{c2} is
 asymptotically equivalent
 to
 $c_0\lambda^{-\frac{d}{r'}+
 2d\left(\frac{1}{2}-\frac{1}{r}\right)}$
 ($c_0>0$),
 as $\lambda\to0^+$, which, combined with the estimates
  \eqref{z02}
 and \eqref{z01}, yields
 $r\geq2$.\par
Similarly, if the inequality \eqref{z1}
holds
 for $t\in(0,\delta)$, one
 must have
 \begin{equation}\label{ul}
 \frac{\left[(1+
 \lambda^{-2})^2+
 t^2\right]^{\frac{d}{2}
 \left(\frac{1}{2}-\frac{1}{r}\right)}}
 {\left[
 \lambda^{-2}(\lambda^{-2}+1)+t^2
 \right]^{d\left(\frac{1}{2}-
 \frac{1}{r}\right)}}\leq
 Ct^\alpha,
 \end{equation}
 for every $t\in (0,\delta)$,
 $\lambda>0$. Choosing
 $t=\lambda^{-1}$, we see
 that,
 when $\lambda\to+\infty$,
 the left-hand side of
 \eqref{ul} is asymptotically
 equivalent to
 $c_1\lambda^{2d\left(\frac{1}{2}-\frac{1}{r}\right)}$ ($c_1>0$),
 and this proves
 \eqref{z2}.\par
 Finally,  choosing $\lambda=1$ and
 letting $t\to+\infty$, we see
 that the expression in \eqref{c2} is asymptotically
 equivalent to  $c_2
 t^{-d\left(\frac{1}{2}-\frac{1}{r}\right)}$,
 which is \eqref{z3}.
\end{proof}
\begin{proposition}{\rm (Sharpness of
\eqref{str}).}\label{prop2}
Assume that for some $1\leq
\alpha,\beta,r\leq\infty$,
$C>0$, the following estimate
holds:
\begin{equation}\label{s1}\|e^{it\Delta}
u_0\|_{W(L^{\alpha},L^\beta)_t
W(\Fur
L^{r^\prime},L^r)_x}\leq C
\|u_0\|_{L^2_x},\quad\forall
u_0\in\mathcal{S}(\R^d).
\end{equation}
Then $r\geq2$. Moreover, if
$q$ is defined by the scaling
relation
$2/q+d/r=d/2$,
then
\begin{equation}\label{s2}
\alpha\leq\frac{q}{2},
\end{equation}
\begin{equation}\label{s3}
\beta\geq q.
\end{equation}
\end{proposition}
\begin{proof}
 We first prove that $r\geq 2$ and \eqref{s3} holds.
We use the family of initial
data $u_0(\lambda x)=
e^{-\pi\lambda^2 |x|^2}$,
$\lambda>0$. The $W(\Fur
L^{r'},L^r)$ norm of the
corresponding solutions
\eqref{c1} is computed in
\eqref{c2}. We use that
expression to estimate from
below the norm in the left
hand side of \eqref{s1}. We
 single out  $g= \chi_{[-1,1]}$, the characteristic
function of the interval
[-1,1], as window function to compute the
$W(L^\alpha,L^\beta)$ norm. Then, for
$y\in\R$,
\begin{equation}\label{s4}
\|\|u(\lambda^2 t,\lambda
\,\cdot)\|_{W(\Fur
L^{r'},L^r)} T_y
g\|_{L^\alpha_t}\gtrsim
\lambda^{-\frac{d}{r'}}(\lambda^{-4}+y^2)^{
-\frac{d}{2}\left(\frac{1}{2}
-\frac{1}{r}\right)},\quad
0<\lambda\leq1.
\end{equation}
Since the left-hand side of
\eqref{s1} is finite, this
estimate for $\lambda=1$
already implies $r\geq 2$.
Now, to compute the $L^\beta$
norm of the expression in
\eqref{s4} we apply the
formula
\begin{equation}\label{s4bis}
\int_{-\infty}^{+\infty}
(\mu+y^2)^\gamma\,dy=c_\gamma\mu^{\frac{1}{2}+\gamma},
\ \mu>0,
\end{equation}
for every $\gamma$, and for a
convenient $c_\gamma\in
(0,\infty]$, independent of
$\mu$ ($c_\gamma<\infty$ if
$\gamma<-\frac{1}{2}$). We
deduce
\begin{equation}\label{s5}
\|u(\lambda^2\,\cdot,\lambda\,\cdot)\|_{
W(L^{\alpha},L^\beta)_t
W(\Fur
L^{r^\prime},L^r)_x}\gtrsim
\lambda^{-\frac{d}{r'}-\frac{2}{\beta}+
2d\left(\frac{1}{2}-\frac{1}{r}\right)}.
\end{equation}
We write the assumption \eqref{s1} for the solution $u(\lambda^2t,\lambda x)$, corresponding to the initial datum $u_0(\lambda x)$. Using the minorization  \eqref{s5},   the trivial
equality
\begin{equation}
\label{s6}
\|u_0(\lambda\,\cdot)\|_{L^2}=
\lambda^{-\frac{d}{2}}\|u_0\|_{L^2},
\end{equation}
and letting $\lambda\to0^+$, we infer
\[
-\frac{d}{r'}-\frac{2}{\beta}+2d
\left(\frac{1}{2}-\frac{1}{r}\right)\geq
-\frac{d}{2},
\]
that is
$\frac{2}{\beta}\leq
\frac{d}{2}-\frac{d}{r}=\frac{2}{q}$, i.e., the estimate
\eqref{s3}.\par
We now prove \eqref{s2}.
Again we use the formula
\eqref{c2} to estimate
\begin{multline}\label{s7}
\|\|u(\lambda^2 t,\lambda
\,\cdot)\|_{W(\Fur
L^{r'},L^r)}
T_y g\|_{L^\alpha_t}\\
\gtrsim
\lambda^{-\frac{d}{r'}}(1+y^2)^{\frac{d}{2}\left(\frac{1}{2}
-\frac{1}{r}\right)
}\left(\int_{y-1}^{y+1}(\lambda^{-2}+t^2)^{-\alpha
d\left(\frac{1}{2}
-\frac{1}{r}\right)}\,dt\right)^{\frac{1}
{\alpha}},\quad \lambda\geq1.
\end{multline}
An application of the
inequality
\[
\int_{-\frac{1}{2}}^\frac{1}{2}(\mu+t^2)^\gamma\,dt\gtrsim\mu^{\frac{1}{2}+\gamma},\quad
0<\mu\leq 1,
\]
allows us to estimate the
expression in \eqref{s7} as
\[
\|\|u(\lambda^2 t,\lambda
\,\cdot)\|_{W(\Fur
L^{r'},L^r)}
T_y g\|_{L^\alpha_t}\\
\gtrsim
\lambda^{-\frac{d}{r'}-\frac{1}{\alpha}+2d\left(\frac{1}{2}
-\frac{1}{r}\right)},\quad
|y|\leq\frac{1}{2},\
\lambda\geq1.
\]
Taking the $L^\beta_y$ norm,
we obtain
\begin{equation}\label{s8}
\|u(\lambda^2\,\cdot,\lambda\,\cdot)\|_{W(L^{\alpha},L^\beta)_t
W(\Fur
L^{r^\prime},L^r)_x}\gtrsim
\lambda^{-\frac{d}{r'}-\frac{1}{\alpha}+
2d\left(\frac{1}{2}-\frac{1}{r}\right)},\quad
\lambda\geq1.
\end{equation}
As above, using the assumption \eqref{s1} for the solution $u(\lambda^2 t,\lambda x)$, corresponding to the initial datum $u_0(\lambda x)$, and the estimates
\eqref{s6} and \eqref{s8}, we
 let
$\lambda\to+\infty$ and infer
\[
-\frac{d}{r'}-\frac{1}{\alpha}+2d
\left(\frac{1}{2}-\frac{1}{r}\right)\leq
-\frac{d}{2},
\]
namely,  $\frac{1}{\alpha}\geq
\frac{d}{2}-\frac{d}{r}=\frac{2}{q}$, that is \eqref{s2}.
\end{proof}

We now discuss the sharpness
of the estimates
\eqref{hom2i}. The following two auxiliary
results are needed.
\begin{lemma}
Let $a>0, b\in\R$. With the
notation in Lemma \ref{lem1},
we have
\begin{equation}
\label{dd1}
\|\phi^{(a+ib)}\|_{W(L^{r_1},L^{r_2})}
\asymp a^{-
\frac{d}{2r_2}}(a+1)^{\frac{d}{2}\left(\frac{1}{r_2}-\frac{1}{r_1}\right)}.
\end{equation}
\end{lemma}
\begin{proof}
We take the Gaussian function
$g(x)=e^{-\pi|x|^2}$ as
window in the definition of
$W(L^{r_1},L^{r_2})$. Then
\[
|\phi^{(a+ib)}(x)T_y
g(x)|=e^{-\pi[(a+1)|x|^2-2xy+|y|^2]},
\]
so that \eqref{lem0} gives
\[
\|\phi^{(a+ib)}T_y
g\|_{L^{r_1}}=[r_1(a+1)]^{-\frac{d}{2r_1}}e^{-\frac{\pi
a}{a+1}|y|^2}.
\]
Now we compute the $L^{r_2}$
norm of this expression
(again by using \eqref{lem0})
and we obtain \eqref{dd1}.
\end{proof}
\begin{lemma} Let $u(t,x)$ be
the solution of \eqref{cp}
with the initial datum
$u_0(x)=e^{-\pi|x|^2}$. Then,
for the solution
$u(\lambda^2t,\lambda x)$
corresponding to the initial
datum $u_0(\lambda x)$,
$\lambda>0$, we have
\begin{equation}\label{dd2}
\|u(\lambda^2
t,\lambda\,\cdot)\|_{W(L^{r_1},L^{r_2})}
\asymp
   \lambda^{-\frac{d}{r_2}}
(1+(t\lambda^2)^2)^{\frac{d}{2}
\left(\frac{1}{r_1}-\frac{1}{2}\right)}
(1+\lambda^2+(
t\lambda^2)^2)^{\frac{d}{2}
\left(\frac{1}{r_2}-\frac{1}{r_1}\right)}.
\end{equation}
\end{lemma}
\begin{proof}
We use the explicit formula
\eqref{c1}, which can be
written as
\[
u(\lambda^2 t,\lambda
x)=(1+4\pi i
t\lambda^2)^{-d/2}\phi^{(a+ib)},
\]
with
\[
a=\frac{\lambda^2}{1+(4\pi
t\lambda^2)^2},\quad
b=-\frac{4\pi
t\lambda^2}{1+(4\pi
t\lambda^2)^2}.
\]
Hence \eqref{dd2} follows
from the estimate \eqref{dd1}.
\end{proof}
\begin{proposition}\label{pd1} Let $u_0(x)=
e^{-\pi|x|^2}$. Then
\begin{equation}\label{dd3bis}
\|e^{it\Delta}
u_0\|_{W(L^{r_1},L^{r_2})}\sim
C
t^{d\left(\frac{1}{r_2}-\frac{1}{2}\right)},\quad
{\it as}\ t\to+\infty,
\end{equation}
for some $C>0$.\par Moreover,
suppose that for some
$t_0\not=0$, $1\leq
r_1,r_2\leq\infty$, $C>0$,
the following estimate holds:
\begin{equation}\label{dd3}
\|e^{it_0\Delta}u_0\|_{W(L^{r_1},L^{r_2})}\leq
C\|u_0\|_{W(L^{r'_1},L^{r'_2})},\quad\forall
u_0\in\mathcal{S}(\R^d).
\end{equation}
Then $r_1\leq r_2$.
\end{proposition}
\begin{proof}
The formula \eqref{dd3bis} is
an immediate consequence of the norm
\eqref{dd2}, with
$\lambda=1$.\par To prove the
remaining part of the
statement, we apply the assumption
\eqref{dd3} to the initial
data $u_0(\lambda
x)=e^{-\pi\lambda^2|x|^2}$,
$\lambda>0$. As a consequence
of the norm expressions computed in \eqref{dd1} and \eqref{dd2},
 we see that, when
$\lambda\to+\infty$, the
left-hand side and right-hand
side of \eqref{dd3} are
asymptotically equivalent to
$c_1
\lambda^{\frac{d}{r_2}-d}$,
and
$c_2\lambda^{-\frac{d}{r'_1}}$
respectively, for some
$c_1,c_2>0$. This implies
$\frac{d}{r_2}-d\leq
-\frac{d}{r'_1}$, namely
$r_1\leq r_2$.
\end{proof}
\begin{proposition}{\rm (Sharpness of
\eqref{hom2i})}.\label{pd2}
Suppose that, for some $1\leq
q_1,q_2,r_1,r_2\leq\infty$,
$C>0$, the following estimate
holds:
\begin{equation}\label{dd4}
\|e^{it\Delta}
u_0\|_{W(L^{q_1},L^{q_2})_t
W(L^{r_1},L^{r_2})_x}\leq
C\|u_0\|_{L^2},\quad \forall
u_0\in\mathcal{S}(\R^d).
\end{equation}
Then we have
\begin{equation}\label{dd5}
\frac{2}{q_1}+\frac{d}{r_1}\geq\frac{d}{2},
\end{equation}
\begin{equation}\label{dd6}
\frac{2}{q_2}+\frac{d}{r_2}\leq\frac{d}{2}.
\end{equation}
In particular, $($because of
$\eqref{dd6})$, $r_2\geq2$.
Finally, it must be
$q_2\geq2$.
\end{proposition}
\begin{proof}
We apply the estimate
\eqref{dd4} to the family of
initial data $u_0(\lambda
x)=e^{-\pi\lambda^2|x|^2}$,
$\lambda>0$. Choose $g=\chi_{[-1,1]}$
 as window function  in the definition
of the space
$W(L^{q_1},L^{q_2})_t$. Then,
for $y\in\mathbb{R}$,
$\lambda\geq1$,
\begin{multline}\label{dd7}
\| u(\lambda^2
t\lambda,\,\cdot)\|_{W(L^{r_1},L^{r_2})}
T_y g(t)\|_{L^{q_1}}\\
\asymp\lambda^\frac{d}{r_2}\left(\int^{y+1}
_{y-1}(\lambda^{-4}+t^2)^{\frac{q_1
d}{2}\left(\frac{1}{r_1}-\frac{1}{2}\right)}
(\lambda^{-2}+t^2)^{\frac{q_1
d}{2}\left(\frac{1}{r_2}-\frac{1}{r_1}\right)}dt\right)^{\frac{1}{q_1}}.
\end{multline}
To estimate this expression
from below, we apply the
easily verified formula
\[
\int_{-\frac{1}{2}}^{\frac{1}{2}}(\mu^2+t^2)
^{\gamma_1}(\mu+t^2)^{\gamma_2}dt\gtrsim\mu^{1+2\gamma_1+\gamma_2},\quad
0<\mu\leq1,
\]
with $\mu=\lambda^{-2}$. Thereby,
\[
\|\,\| u(\lambda^2
t, \lambda\,\cdot)\|_{W(L^{r_1},L^{r_2})}
T_y g(t)\|_{L^{q_1}}\gtrsim
\lambda^{-\frac{2}{q_1}-\frac{d}{r_1}}\quad
{\rm for}\
|y|\leq\frac{1}{2},\
\lambda\geq1.
\]
Taking  the $L^{q_2}$ norm of this expression yields
\[
\| u(\lambda^2
\,\cdot,\lambda\,\cdot)\|_{{W(L^{r_1},L^{r_2})}_tW(L^{r_1},L^{r_2})_x}
\gtrsim
\lambda^{-\frac{2}{q_1}-\frac{d}{r_1}},
\quad\lambda\geq1.
\]
On the other hand, the
right-hand side of
\eqref{dd4} is equal to
$\lambda^{-\frac{d}{2}}$,
therefore letting
$\lambda\to+\infty$, we obtain the index relation
\eqref{dd5}.\par We now prove
\eqref{dd6}. If
$0<\lambda\leq1$ it follows
from the norm estimate \eqref{dd2} that
\[
\|u(\lambda^2
t,\lambda\,\cdot)\|_{W(L^{r_1},L^{r_2})}\asymp\lambda^{\frac{d}{r_2}-d}(\lambda^{-4}+t^2)
^{\frac{d}{2}\left(\frac{1}{r_2}-
\frac{1}{2}\right)}.
\]
Hence, we have
\[
\|\|u(\lambda^2
t,\lambda\,\cdot)\|_{W(L^{r_1},L^{r_2})}T_y
g\|_{L^{q_1}}\asymp\lambda^{\frac{d}{r_2}-d}
(\lambda^{-4}+y^2)
^{\frac{d}{2}\left(\frac{1}{r_2}-\frac{1}{2}
\right)},\quad
0<\lambda\leq1.
\]
 Taking
the $L^{q_2}$ norm of
this expression and
applying the formula \eqref{s4bis},
with $\mu=\lambda^{-4}$, we obtain the minorization
\[
\|u(\lambda^2
\,\cdot,\lambda\,\cdot)\|
_{W(L^{q_1},L^{q_2})_tW(L^{r_1},L^{r_2})_x}\gtrsim\lambda^{-\frac{2}{q_2}-\frac{d}{r_2}},\quad
0<\lambda\leq1.
\]
Since the right-hand
side of \eqref{dd4} is equal
to $\lambda^{-\frac{d}{2}}$, the estimate \eqref{dd6} follows by
letting $\lambda\to0^+$.\par We
are left to prove the condition $q_2\geq
2$. We follow the pattern
outlined in \cite[Exercise 2.42]{tao}. Precisely, we
take as initial datum
\[
u_0(x)=\sum_{j=1}^N
e^{-it_j\Delta} f,
\]
where $f$ is a fixed test
function, normalized so that the function $v(t,x)$, defined by
$v(t,x)=\left(e^{it\Delta}
f\right)(x)$, satisfies
\begin{equation}\label{zi0}
\|v\|_{W(L^{q_1},L^{q_2})_t
W(L^{r_1},L^{r_2})_x}=1.
\end{equation}
Here $t_1,\ldots,t_N$ are
widely separated times that
will be chosen later on.
Notice that the corresponding
solution will be
$$ u(t,x)=\left(e^{it\Delta}
u_0\right)(x)=\sum_{j=1}^N
v(t-t_j,x).$$
 We claim
that
\begin{equation}\label{zi1}
\|u_0\|_{L^2}=\|\sum_{j=1}^N
e^{-it_j\Delta}
f\|_{L^2}\leq(N+1)^{\frac{1}{2}}\|f\|_{L^2},
\end{equation}
if $t_1,\ldots,t_N$ are
suitable separated.\\
 Indeed,
 \begin{align}
\|\sum_{j=1}^N
e^{-it_j\Delta}
f\|_{L^2}^2&=\sum_{j=1}^N
\|e^{-it_j\Delta}
f\|_{L^2}^2+\sum_{j\not=k}\langle
e^{i(t_j-t_k)\Delta}
f,f\rangle\nonumber\\
&\leq
N\|f\|_{L^2}^2+C\sum_{j\not=k}
|t_j-t_k|^{-\frac{d}{2}}\|f\|_{L^1}^2,\nonumber
\end{align}
where we used Cauchy-Schwarz
inequality and the classical
dispersive estimate
$\|e^{it\Delta}f\|_{L^\infty}\leq
C|t|^{-d/2}\|f\|_{L^1}$.
Hence \eqref{zi1} follows if
$$|t_j-t_k|^{-d/2}\leq[C(N^2-N)\|f\|_{L^1}^2]^{-1}
\|f\|_{L^2}^2.$$\par We now
estimate from below the
left-hand side of
\eqref{dd4}. To this end, let
$\tilde{v}(t,x)=v(t,x)\chi_R(t)$,
where $\chi_R(t)$ is the
characteristic function of
the interval $[-R,R]$.
Moreover, we assume
$q_2<\infty$ and  choose $R$
large enough so that
\begin{equation}\label{zi2}
\|v-\tilde{v}\|_{W(L^{q_1},L^{q_2})_t
W(L^{r_1},L^{r_2})_x}\leq\frac{1}{N}.
\end{equation}
 We claim that, if
 $|t_j-t_k|\geq 2R+2$, for every $ j\not=k$, then
 \begin{eqnarray}
\|u(t,x)\|_{W(L^{q_1},L^{q_2})_t
W(L^{r_1},L^{r_2})_x}&=&\|\sum_{j=1}^N
v(t-t_j,x)\|_{W(L^{q_1},L^{q_2})_t
W(L^{r_1},L^{r_2})_x}\nonumber\\
&\geq&
N^{\frac{1}{q_2}}\!\left(1-\frac{1}{N}\right)-1.\label{zi3}
 \end{eqnarray}
 This,  together with the assumption \eqref{dd4} and
 the $L^2$-estimate
 of the initial  datum  \eqref{zi1},
 gives the condition
 $q_2\geq2$, for $N$ large
 enough.\par
 In order to prove the minorization
 \eqref{zi3}, observe that, by the assumption
 \eqref{zi2},
 \[
 \|\sum_{j=1}^N
 v(t-t_j,x)-\sum_{j=1}^N
 \tilde{v}(t-t_j,x)\|_{W(L^{q_1},L^{q_2})_t
W(L^{r_1},L^{r_2})_x}\leq 1.
\]
Hence it suffices to prove
\begin{equation}\label{zi4}
\|\sum_{j=1}^N
\tilde{v}(t-t_j,x)\|_{W(L^{q_1},L^{q_2})_t
W(L^{r_1},L^{r_2})_x}\geq
N^{\frac{1}{q_2}}\left(1-\frac{1}{N}\right).
 \end{equation}
Now,
\[
\|\sum_{j=1}^N
\tilde{v}(t-t_j,\cdot)\|_{
W(L^{r_1},L^{r_2})}=\sum_{j=1}^N
\|\tilde{v}(t-t_j,\cdot)\|_{
W(L^{r_1},L^{r_2})},
\]
since, for every fixed $t$,
there is at most one function
in the sum which is not
identically zero. Hence, upon
setting
$h_j(t):=\|\tilde{v}(t-t_j,\cdot)\|_{
W(L^{r_1},L^{r_2})}$,  we
have
\[
\|\sum_{j=1}^N
\tilde{v}(t-t_j,\cdot)\|_{W(L^{q_1},L^{q_2})_t
W(L^{r_1},L^{r_2})_x}=\|\|\sum_{j=1}^N
h_j T_y
g\|_{L^{q_1}}\|_{L^{q_2}_y}.
\]
Choosing the window function $g$
  supported in
$[0,1]$, since the $h_j$'s
are supported in intervals
separated by a distance
$\geq2$, we see that the last
expression is equal to
\[
\|\sum_{j=1}^N \|h_j T_y
g\|_{L^{q_1}}\|_{L^{q_2}_y}.
\]
In turn, since the functions
$y\longmapsto\|h_j T_y g\|_{L^{q_1}}$,
$j=1,\ldots,N$, have disjoint supports,
the norm above can be written as
\[
\left(\sum_{j=1}^N \|\|h_j T_y
g\|_{L^{q_1}}\|_{L^{q_2}_y}^{q_2}
\right)^{\frac{1}{q_2}}=N^{\frac{1}{q_2}}
\|\tilde{v}\|_{W(L^{q_1},L^{q_2})_t
W(L^{r_1},L^{r_2})_x}.
\]
 Hence, the minorization  \eqref{zi4} follows
from the assumptions \eqref{zi0} and \eqref{zi2}.
\end{proof}

\section{An application to Schr\"odinger equations with
 time-dependent potentials}
In this section we prove the
wellposedness in $L^2$ of
following Cauchy problem in
any dimension $d\geq 1$:
\begin{equation}\label{cpP}
\begin{cases}
i\partial_t u+\Delta u=V(t,x)u,\quad t\in [0,T]=I_T, \,\,x\in\rd,\\
u(0,x)=u_0(x),
\end{cases}
\end{equation}
for the class of potentials
\begin{equation}\label{pot3}
V\in L^\alpha(I_T;W(\cF
L^{p^\prime},
L^p)_x),\quad\frac1{\a}+\frac{d}{p}\leq1,\
1\leq\alpha<\infty,\ {d}<
p\leq\infty.
\end{equation}
Precisely, we generalize
\cite[Theorem 6.1]{cordero}
by treating the one
dimensional case as well and
allowing the potential to
belong to Wiener amalgam
space with respect to $x$,
rather than simply $L^p$
spaces.
\par To this end, let us
prove directly a simple
point-wise multiplication
property of Wiener amalgam
spaces (it is a special case
of the \modsp \, property
\cite[Proposition
2.4]{CG02}).

\begin{lemma}
 Let  $1\leq p,q,r\leq\infty$. If
\begin{equation}\label{ind}\frac1p+\frac1q=\frac1{r^\prime},\end{equation}
then
\begin{equation}\label{mconvm}
W(\cF L^{p^\prime},L^p)(\Ren)
\cdot W (\cF
L^{q\prime},L^q)(\Ren)\subset
W(\cF L^r,L^{r^\prime})(\Ren)
\end{equation}
with  norm inequality  $\| f
h \|_{ W(\cF
L^r,L^{r^\prime})}\lesssim
\|f\|_{W(\cF
L^{p^\prime},L^p)}\|h\|_{W
(\cF L^{q^\prime},L^q)}$.
\end{lemma}
\begin{proof}
We measure the Wiener amalgam
norm with respect to the
window $g(x) =g_1(x) g_2(x)$,
$g_1,g_2\in
\cC_0^\infty(\rd)$, with
$\|g\|_2=\|g_1\|_2=\|g_2\|_2=1$.
(Different windows  yield
equivalent norms for the
Wiener amalgam spaces).

Using
$\hat{g}=\hat{g_1}\ast\hat{g_2}$,\,
$\widehat{T_xg}=M_{-x}\hat{g}$
and
$$M_{-x} \hat{g}=M_{-x} \hat{g_1}\ast M_{-x} \hat{g_2},
$$
we can write
\begin{eqnarray*}
\|f h\|_{W(\cF L^r,L^{r^\prime})}&\asymp & \| \|\hat{f}\ast\hat{h}\ast \widehat{T_x g}\|_{L^r}\|_{L^{r^\prime}}\\
&=& \| \|(\hat{f}\ast M_{-x}\hat{g_1})\ast(\hat{h}\ast (M_{-x}\hat{ g_2})\|_{L^r}\|_{L^{r^\prime}}\\
&=& \| \|\widehat{fT_{x}{g_1}}\ast\widehat{hT_{x}{g_2}}\|_{L^r}\|_{L^{r^\prime}}\\
&\lesssim & \| \|fT_{x}{g_1}\|_{\cF L^{p^\prime}}\|hT_{x}{g_2}\|_{\cF L^{q^\prime}}\|_{L^{r^\prime}}\\
&\lesssim & \| \|fT_{x}{g_1}\|_{\cF
L^{p^\prime}}\|_{L^p}\|\|hT_{x}{g_2}\|_{\cF L^{q^\prime}}\|_{L^{q}}\\
&=& \|f\|_{W(\cF
L^{p\prime},L^p)}\|h\|_{W
(\cF L^{q\prime},L^q)},
\end{eqnarray*}
where the  former inequality
is the consequence of Young's
Inequality with
$1/{p^\prime}+1/{q^\prime}=1+1/r$,
which follows from the
assumption \eqref{ind}, and
the latter is H\"older's
inequality with index
relation \eqref{ind}.
\end{proof}
\vskip0.3truecm We have now
the instruments to prove the
following result.
\begin{theorem}\label{tepot}
Consider the class of
potentials \eqref{pot3}.
Then, for all $(q,r)$ such that $2/q+d/r=d/2$, $q>4,r\geq2$, the Cauchy problem \eqref{cpP} has a
 unique solution \\
(i) $u\in
\mathcal{C}(I_T;L^2(\bR))
\cap L^{q/2}(I_T; W(\cF
L^{r^\prime}, L^r))$, if
$d=1$;
\\
(ii) $u\in
\mathcal{C}(I_T;\lrd)\cap
L^{q/2}(I_T; W(\cF
L^{r^\prime}, L^r))\cap
L^{2}(I_T; W(\cF
L^{2d/(d+1),2},
L^{2d/(d-1)}))$, if $d>1$.
\end{theorem}
\begin{proof}
It is enough to prove the
case $d=1$. Indeed, for
$d\geq 2$, condition
\eqref{pot3} implies $p>2$,
so that  $\cF L^{p^\prime}
\hookrightarrow L^p$ and the
inclusion relations of Wiener
amalgam spaces yield $W(\cF
L^{p^\prime},L^p)\hookrightarrow
W( L^{p},L^p)=L^p$. Hence our
class of potentials is a
subclass of those of
\cite[Theorem 6.1]{cordero},
for which the quoted theorem
provides the desired
result.\par

We now turn to the case
$d=1$. The proof follows the
ones of \cite[Theorem 1.1,
Remark 1.3]{Dancona05} and
\cite[Theorem 6.1]{cordero}
(see also \cite{Yajima87}).
\par First of all, since the
interval $I_T$ is bounded, by
H\"older's inequality it
suffices to assume
$1/\alpha+d/p=1$.\par
 We
choose a small time interval
$J=[0,\delta]$  and set, for
$q\geq2$, $q\not=4$, $r\geq
1$,
\[
Z_{q/2,r}=L^{q/2}(J;W(\cF
L^{r^\prime}, L^r)_x).
\]
Now,  fix  an admissible pair
$(q_0,r_0)$ with $r_0$
arbitrarily large (hence
$(1/q_0,1/r_0)$ is
arbitrarily close to
$(1/4,0)$) and set
$Z=\mathcal{C}(J;L^2)\cap
Z_{q_0/2,r_0}$, with the norm
$\|v\|_Z=\max\{\|v\|_{\mathcal{C}(J;L^2)},\|v\|_{Z_{q_0/2,r_0}}\}$.
We have $Z\subset Z_{q/2,r}$
for all  admissible pairs
$(q,r)$
 obtained
by interpolation between
$(\infty,2)$ and $(q_0,r_0)$.
Hence, by the arbitrary of
$(q_0,r_0)$ it suffices to
prove that $\Phi$ defines a
contraction in $Z$.\par
 Consider now
the integral formulation of
the Cauchy problem, namely
$u=\Phi(v)$, where
\[
\Phi(v)=e^{it\Delta}u_0+\int_0^t
e^{i(t-s)\Delta}V(s)v(s)\,ds.
\]
By the homogeneous and
retarded Strichartz estimates
in \cite[Theorems 1.1,
1.2]{cordero} the following
inequalities hold:
\begin{equation}\label{nan}
\|\Phi(v)\|_{Z_{q/2,r}}\leq
C_0
\|u_0\|_{L^2}+C_0\|Vv\|_{Z_{(\tilde{q}/2)',\tilde{r}'}},
\end{equation}
for all admissible pairs
$(q,r)$ and
$(\tilde{q},\tilde{r})$,
$q>4,\tilde{q}>4$.\par
Consider now the case $1\leq
\alpha< 2$. We choose
$((\tilde{q}/2)',\tilde{r})=(\alpha,2p/(p+2))$.
Since $v\in L^\infty(J; L^2)$
applying \eqref{mconvm} for
$q=2$ we get
$$\|Vv\|_{W(\cF L^{\tilde{r}},L^{\tilde{r}^\prime})}\lesssim \|V\|_{W(\cF
L^{p^\prime},L^{p})}\|v\|_{L^2},$$
whereas H\"older's Inequality
in the time-variable gives
$$\|Vv\|_{Z_{(\tilde{q}/2)',\tilde{r}'}}\lesssim \|V\|_{L^{\alpha}(J;W(\cF
L^{p^\prime},L^{p}))}\|v\|_{L^\infty(J;L^2)}.
$$
The estimate \eqref{nan} then
becomes
\begin{equation}\label{nan2}
\|\Phi(v)\|_{Z_{q/2,r}}\leq
C_0
\|u_0\|_{L^2}+C_0\|V\|_{L^{\alpha}(J;W(\cF
L^{p^\prime},L^{p}))}\|v\|_{L^\infty(J;L^2)}.
\end{equation}
By taking $(q,r)=(\infty,2)$
or    $(q,r)=(q_0,r_0)$ one
deduces that $\Phi:Z\to Z$
(the fact that $\Phi(u)$ is
{\it continuous} in $t$ when
valued in $L^2_x$ follows
from a classical limiting
argument \cite[Theorem 1.1,
Remark 1.3]{Dancona05}).
 Also,  if $J$ is small enough,
$C_0\|V\|_{L^\alpha_t
L^p_x}<1/2$, and $\Phi$ is a
contraction. This gives a
unique solution in $J$. By
iterating this argument a
finite number of times one
obtains a solution in
$[0,T]$.

 \par The case $2\leq \alpha<\infty$ is similar.
We again consider the
inequality \eqref{nan} for
all admissible pairs $(q,r)$
and $(\tilde{q},\tilde{r})$,
$q>4$, $\tilde{q}>4$. Since
$\alpha \geq 2$ we can find
an admissible pair
$(\tilde{q},\tilde{r})$,
$\tilde{q}>4$, such that
\begin{eqnarray}
\frac{1}{{(\tilde{q}/2)}^\prime}&=&\frac{1}{(q_0/{2})}+\frac1{\alpha}\label{ind1}\\
\frac{1}{\tilde{r}^\prime}&=&\frac{1}{r_0}+\frac1{p}\label{ind2}\\
\end{eqnarray}
 Using the Wiener point-wise property \eqref{mconvm} with index relation \eqref{ind2}  we have
$$
\|Vv\|_{W(\cF
L^{\tilde{r}},L^{\tilde{r}^\prime})}\lesssim
\| V\|_{W(\cF
L^{p^\prime},L^{p})} \|v \|
_{W(\cF
L^{r_0^\prime},L^{r_0})}.
$$
Finally, H\"older's
Inequality  with index
relation \eqref{ind1}  gives
$$
\|Vv\|_{L^{(\tilde{q}/2)'}(J;W(\cF
L^{\tilde{r}},L^{\tilde{r}^\prime}))}\lesssim
\| V\|_{L^\alpha(J; W(\cF
L^{p^\prime},L^{p}))} \|v
\|_{L^{q_0/2}(J;W(\cF
L^{r_0^\prime},L^{r_0}))}.
$$
The final part of the proof
is analogous to the previous
case.
\end{proof}

\begin{Remark}\label{espot} The most interesting case in Theorem \ref{tepot} is when $d=1$. Indeed, choosing $p>2$, so that $1<p'<2$, we have the embedding $H_s\hookrightarrow \cF L^{p'}$, for $s>1/{p'}-1/2$ (see, e.g., \cite[Theorem
7.9.3]{Hormander1}). Whence, examples of potentials $V(t,x)$ satisfying \eqref{pot3} are given by tempered distributions locally in $H_s$ as above in the $x$-variable, conveniently localized in $x$, and belonging to $L^\alpha_t$, $\alpha\geq p'$, with respect to the $t$-variable.
\end{Remark}

\appendix
\section{Some
results in real interpolation
theory} Here we collected
some results in real
interpolation theory which
are used in the proof of the
Strichartz estimates.\par
 Let
$(X,\cB,\mu)$ a measure
space, where $X$ is a set,
$\cB$ a $\sigma$-algebra and
$\mu$ a positive
$\sigma$-finite measure. If
$A$ is a Banach space, $1\leq
p\leq \infty$, then we shall
write $L^p(A)$ for the usual
vector-valued $L^p$ spaces in
the sense of the Bochner
integral. The first result is
a generalization of
\cite[Proposition
2.3]{cordero}. The proof uses
 arguments similar to those in \cite[Section
1.18.4]{triebel}.
\begin{proposition}\label{interlp}
Let $\{A_0,A_1\}$ be an
interpolation couple. For
every $1\leq p_0,p_1<\infty$,
$0<\theta<1$,
$1/p=(1-\theta)/p_0+\theta/p_1$
and $p\leq q$ we have
\begin{equation}\label{inter1}
L^p\left((A_0,A_1)_{\theta,q}\right)\hookrightarrow
(L^{p_0}(A_0),L^{p_1}(A_1))_{\theta,q}.
\end{equation}
\end{proposition}
\begin{proof}
We use the fact that the
bounded functions with values
in $A_0\cap A_1$, vanishing
outside a set with finite
measure, are dense in
$L^{p_0}(A_0)\cap
L^{p_1}(A_1)$ and may be
approximated in
$L^{p_0}(A_0)\cap
L^{p_1}(A_1)$ by simple
functions:
$$c(x)=\sum_{j=1}^N c^{j}\chi_{F_j}(x),\quad \mu(F_j)<\infty,\quad
c^{j}\in A_0\cap A_1, \quad
F_j\cap F_k=\emptyset,
\,\,\mbox{if}\,\,j\not= k.
$$

Set $\eta=p\theta/p_1$ and
$q=pr$, $r\geq1$. If $c(x)$
is such a function, it
follows from Theorem 1.4.2 of
\cite{triebel}, page 29, that
\begin{align}
\|c\|_{(L^{p_0}(A_0),L^{p_1}(A_1))_{\theta,q}}^p &\asymp\|t^{-\eta}\inf_{\substack{a+b=c\\
a\in L^{p_0}(A_0),\ b\in
L^{p_1}(A_1)}}
\|a\|_{L^{p_0}(A_0)}^{p_0}+t\|b\|_{L^{p_1}(A_1)}^{p_1}\|_{L^r(\mathbb{R}_+,\frac{dt}{t})}\nonumber\\
&=\|t^{-\eta}
\inf_{\substack{a+b=c\\ a\in
L^{p_0}(A_0),\ b\in
L^{p_1}(A_1)}}
\int_{X} (\|a(x)\|^{p_0}_{A_0}+t\|b(x)\|^{p_1}_{A_1})\|_{{L^{r}(\mathbb{R}_+,\frac{dt}{t})}}\nonumber \\
 &=\|t^{-\eta}\int_{X}\inf_{\substack{a(x)+b(x)=c(x)\\
a(x)\in A_0 ,\ b(x)\in A_1}}
(\|a(x)\|^{p_0}_{A_0}+t\|b(x)\|^{p_1}_{A_1})\|_{{L^{r}(\mathbb{R}_+,\frac{dt}{t})}}\label{ugual}.
\end{align}
In order to justify the last
equality we observe that the
inequality ``$\geq$'' is
obvious. As far as the
converse inequality $``\leq"$
concerns, we take advantage
of the fact that $c(x)$ is a
simple function, so that one
can find minimizing sequences
$a_n(x)$, $b_n(x)$ given by
 simple functions which are constants where $c(x)$ is, and zero
  where $c(x)=0$.\par
Finally, by Minkowski's
integral inequality,
\begin{align}
\|c\|_{(L^{p_0}(A_0),L^{p_1}(A_1))_{\theta,q}}^p&\lesssim\int_{X}\|t^{-\eta}
\inf_{\substack{a(x)+b(x)=c(x)\\
a(x)\in A_0 ,\ b(x)\in A_1}}(\|a(x)\|^{p_0}_{A_0}+t\|b(x)\|^{p_1}_{A_1})\|_{{L^{r}(\mathbb{R}_+,\frac{dt}{t})}}\nonumber\\
&\asymp\|c\|^p_{L^p\left((A_0,A_1)_{\theta,q}\right)}.
\end{align}
\end{proof}

We also record
\cite[Proposition
2.5]{cordero}.
\begin{proposition}\label{inter9}
Given two local components
$B_0,B_1$, for every $1\leq
p_0,p_1<\infty$,
$0<\theta<1$,
$1/p=(1-\theta)/p_0+\theta/p_1$,
and $p\leq q$ we have
\[
W\left((B_0,B_1)_{\theta,q},L^p\right)\hookrightarrow\left(W(B_0,L^{p_0}),W(B_1,L^{p_1})\right)_{\theta,q}.
\]
\end{proposition}
We finally recall
\cite[Section
3.13.5(b)]{bergh-lofstrom}.
\begin{lemma}\label{interp}
If $A_i, B_i, C_i$, $i=0,1$,
are Banach spaces and T is a
bilinear operator bounded
from
\begin{equation*}
    T:A_0\times B_0\rightarrow C_0
\end{equation*}
\begin{equation*}
    T:A_0\times B_1\rightarrow C_1
\end{equation*}\begin{equation*}
    T:A_1\times B_0\rightarrow C_1,
\end{equation*}
then, if $0<\theta_i<\theta<1
$, $i=0,1$,
$\theta=\theta_0+\theta_1$,
 one has
$$T: (A_0,A_1)_{\theta_0,2}\times
(B_0,B_1)_{\theta_1,2}\rightarrow
(C_0,C_1)_{\theta,1}.
$$
\end{lemma}
\section*{Acknowledgements}
The authors would like to
thank Professors Piero
D'Ancona and Luigi Rodino for
fruitful conversations and
comments.


\begin{thebibliography}{10}
\bibitem{benyi} A. B\'enyi, K. Gr\"ochenig,
 K.A. Okoudjou and L.G. Rogers.
 Unimodular Fourier multipliers for
 modulation spaces. {\it J. Funct. Anal.},
246(2): 366-384, 2007.
\bibitem{benyi2} A. B\'enyi and K.A. Okoudjou. Time-frequency estimates for pseudodifferential operators. {\it Contemporary Math.}, Amer. Math. Soc., 428:13--22, 2007.
\bibitem{bergh-lofstrom}
J.~Bergh and J.~L{\"o}fstr{\"o}m.
\newblock {\em Interpolation {S}paces. {A}n {I}ntroduction}.
\newblock Springer-Verlag, Berlin, 1976.
\newblock Grundlehren der Mathematischen Wissenschaften, No. 223.
\bibitem{CG02} E.~Cordero and
K.~Gr\"ochenig.
\newblock {T}ime-frequency analysis
of {L}ocalization operators.
\newblock {\em J. Funct. Anal.},
205(1):107--131, 2003.
\bibitem{cordero}
E.~Cordero and F. Nicola.
\newblock {S}trichartz
estimates in Wiener amalgam spaces for
the Schr\"odinger equation.
\newblock {\em Math. Nachr.},
to appear (available at
ArXiv:math.AP/0610229).
\bibitem{cordero2}
E.~Cordero and F. Nicola.
\newblock {Metaplectic representation
on
 Wiener amalgam spaces and applications
to the Schr\"odinger
equation}. Preprint May 2007,
(available at
ArXiv:0705.0920v1).
\bibitem{Dancona05} P. D'Ancona, V. Pierfelice  and N. Visciglia.
             Some remarks on the {S}chr\"odinger equation with a potential
              in {$L\sp r\sb tL\sp s\sb x$}. {\it Math. Ann.}, 333(2):271--290, 2005.
\bibitem{feichtinger80} H.~G.
Feichtinger.
\newblock Banach convolution algebras of {W}iener's type,
\newblock In {\em Proc. Conf. ``Function, Series, Operators",
 Budapest August 1980}, Colloq. Math. Soc. J\'anos Bolyai, 35,  509--524, North-Holland, Amsterdam,
 1983.
\bibitem{feichtinger83}
H.~G. Feichtinger.
\newblock Banach spaces of distributions of {W}iener's type and interpolation.
\newblock In {\em Proc. Conf. Functional Analysis and Approximation,
 Oberwolfach August 1980},  Internat. Ser. Numer. Math., 69:153--165. Birkh\"auser, Boston, 1981.
\bibitem{feichtinger90}
H.~G. Feichtinger.
\newblock Generalized amalgams, with applications to {F}ourier transform.
\newblock {\em Canad. J. Math.}, 42(3):395--409, 1990.
\bibitem{fournier-stewart85}
J.~J.~F. Fournier and J.~Stewart.
\newblock Amalgams of ${L}\sp p$ and $l\sp q$.
\newblock {\em Bull. Amer. Math. Soc. (N.S.)}, 13(1):1--21, 1985.
\bibitem{Fei98}
H.~G. Feichtinger and G.~Zimmermann.
\newblock {A} {B}anach space of test functions for {G}abor analysis.
\newblock In H.~G. Feichtinger and T. Strohmer, editors, {\em Gabor Analysis and Algorithms. Theory and Applications}, Applied and Numerical Harmonic Analysis, 123--170. Birkh\"auser, Boston, 1998.
\bibitem{folland}
G.~B. Folland.
\newblock {\em Harmonic Analysis in Phase Space}.
\newblock Princeton Univ. Press, Princeton, NJ, 1989.
\bibitem{GinibreVelo92}
J.~Ginibre and G.~Velo.
\newblock Smoothing properties and retarded estimates for some dispersive
  evolution equations.
\newblock {\em Comm. Math. Phys.}, 144(1):163--188, 1992.
\bibitem{grochenig} K. Gr\"ochenig. {\it Foundation of Time-Frequency Analysis}. Birkh\"auser, Boston MA, 2001.
\bibitem{Heil03}
C. Heil.
\newblock {A}n introduction to weighted {W}iener amalgams.
\newblock In M. Krishna, R. Radha and S. Thangavelu, editors, {\em
Wavelets and their
Applications}, 183--216.
Allied Publishers Private
Limited, 2003.
\bibitem{Hormander1}
L.~H{\"o}rmander.
\newblock {\em The analysis of linear partial differential operators. {I}},
  Second Edition.
\newblock Springer-Verlag, Berlin, 1990.
\bibitem{Kato70}
T.~Kato.
\newblock Linear evolution equations of ``hyperbolic'' type.
\newblock {\em J. Fac. Sci. Univ. Tokyo Sect. I}, 17:241--258, 1970.
\bibitem{keel}
M. Keel and T. Tao.
\newblock { Endpoint Strichartz estimates}.
\newblock {\it Amer. J. Math.}, 120:955--980, 1998.
\bibitem{nicola} F. Nicola. Remarks on
dispersive estimates and
curvature. {\it Commun. Pure
Appl. Anal.}, 6:203--212,
2007.
\bibitem{stein}
E. M. Stein.
\newblock {\it Singular integrals and differentiability properties of functions}.
\newblock Princeton University Press, Princeton, 1970.
\bibitem{stein93}
E. M. Stein.
\newblock {\it Harmonic analysis}.
\newblock Princeton University Press, Priceton,1993.
\bibitem{steinweiss}
E. M. Stein and G. Weiss.
\newblock {\it Introduction to Fourier Analysis on Euclidean spaces}.
\newblock Princeton University Press, 1971.
\bibitem{strichartz} R. S.
Strichartz. Restriction of
Fourier transform to
quadratic surfaces and decay
of solutions of wave
equations. {\it Duke Math.
J.}, 44:705-774, 1977.
\bibitem{sugimototomita}
M. Sugimoto and N. Tomita.
\newblock {{T}he  dilation property of modulation spaces and their inclusion relation with Besov spaces}.
\newblock {\em J. Funct. Anal.}, to appear.
\bibitem{tao3} T. Tao.
\newblock {Low regularity
semilinear wave equations}. {\it Comm.
Partial Differential Equations},
24:599-629, 1999.
\bibitem{tao2}
T. Tao.
\newblock {Spherically averaged endpoint  Strichartz estimates for the two-dimensional Schr\"odinger equation}. {\it Comm.\ Partial differential Equations}, 25:1471--1485, 2000.
\bibitem{tao}
T. Tao.
\newblock {\it Nonlinear Dispersive Equations: Local and Global Analysis}.
\newblock CBMS Regional Conference Series in Mathematics, Amer. Math. Soc., 2006.
\bibitem{Toft04}
J.~Toft.
\newblock Continuity properties for modulation spaces, with applications to
  pseudo-differential calculus. {I}.
\newblock {\em J. Funct. Anal.}, 207(2):399--429, 2004.
\bibitem{triebel}
H. Triebel.
\newblock {\it Interpolation theory, function spaces, differential operators}.
\newblock North-Holland, 1978.
\bibitem{baoxiang3} B. Wang, C. Huang.
Frequency-uniform
decomposition method for the
generalized BO, KdV and NLS
equations. {\it J.
Differential Equations}, to
appear.
\bibitem{baoxiang2} B. Wang, H. Hudzik.
The global
 Cauchy problem for the NLS
  and NLKG with small rough data. {\it J.
Differential Equations},
232:36--73, 2007.
\bibitem{baoxiang} B. Wang, L. Zhao and
B. Guo. Isometric
decomposition operators,
function spaces
$E_{p,q}^\lambda$ and
applications to nonlinear
evolution equations. {\it J.
Funct. Anal.}, 233(1):1--39,
2006.
\bibitem{Yajima87}
K.~Yajima.
\newblock Existence of solutions for {S}chr\"odinger evolution equations.
\newblock {\em Comm. Math. Phys.}, 110(3):415--426, 1987.
\end{thebibliography}
\end{document}